\input ./style/arxiv-general.cfg
\documentclass[MSNbibl,number,citesort,seceqn,dvips]{arxbj}
\makeatletter
   \@ifpackageloaded{graphicx}{}{\usepackage{graphicx}}
\makeatother
\usepackage{upgreek}

\volume{22}
\issue{1}
\pubyear{2016}
\firstpage{615}
\lastpage{651}
\doi{10.3150/14-BEJ670} 
\docsubty{FLA}

\makeatletter
\newcommand{\mathsfTheta}{\ooalign{\textsf{O}\cr\hidewidth\raisebox{1pt}{\textbf{\textup{-}}}\hidewidth\cr}}
\newcommand{\KL}{\mathrm{KL}}
\newcommand{\rrVert}{\Vert}
\newcommand{\llVert}{\Vert}
\newcommand{\rrvert}{\vert}
\newcommand{\llvert}{\vert}
\newcommand{\argmax}{\mathop{\arg\max}}
\newcommand{\argmin}{\mathop{\arg\min}}
\newcommand{\I}{\mathrm{I}}
\newcommand{\II}{I}
\newcommand{\Var}{\operatorname{Var}}
\newcommand{\Cov}{\operatorname{Cov}}
\newcommand{\tr}{\operatorname{tr}}
\newcommand{\diam}{\operatorname{diam}}
\newcommand{\asto}{\mathop{\longrightarrow}\limits^{\mathrm{a.s.}}}
\newcommand{\astoi}{\mathop{\longrightarrow}\limits^{\mathit{a.s.}}}
\newcommand{\prhoto}{\mathop{\leadsto}\limits^{\mathrm{P}}}
\newcommand{\pto}{\mathop{\to}\limits^{\mathrm{P}}}
\renewcommand{\Pr}{\mathsf{P}}

\newtheorem{theorem}{Theorem}[section]
\newtheorem{lemma}{Lemma}[section]
\newtheorem{corollary}{Corollary}[section]
\newproclaim{definition}{Definition}[section]
\newremark{remark}{Remark}[section]
\newremark{example}{Example}[section]
\makeatother

\begin{document}
\begin{frontmatter}

\title{Asymptotic optimality of myopic information-based strategies for Bayesian adaptive estimation}
\runtitle{Asymptotic optimality of myopic strategies}

\begin{aug}
\author{\inits{J.V.}\fnms{Janne V.}~\snm{Kujala}\corref{}\ead[label=e1]{jvk@iki.fi}}
\address{Department of Mathematical Information Technology,
University of Jyv{\"a}skyl{\"a}, P.O. Box 35, FI-40014 Jyv{\"a}skyl{\"
a}, Finland. \printead{e1}}
\end{aug}

\received{\smonth{3} \syear{2012}}
\revised{\smonth{5} \syear{2014}}

%
\begin{abstract}
This paper presents a general asymptotic theory of sequential
Bayesian estimation giving results for the strongest, almost sure
convergence. We show that under certain smoothness conditions on
the probability model, the greedy information gain maximization
algorithm for adaptive Bayesian estimation is asymptotically optimal
in the sense that the determinant of the posterior covariance in a
certain neighborhood of the true parameter value is asymptotically
minimal. Using this result, we also obtain an asymptotic expression
for the posterior entropy based on a novel definition of almost sure
convergence on ``most trials'' (meaning that the convergence holds
on a fraction of trials that converges to one). Then, we extend the
results to a recently published framework, which generalizes the
usual adaptive estimation setting by allowing different trial
placements to be associated with different, random costs of
observation. For this setting, the author has proposed the
heuristic of maximizing the expected information gain divided by the
expected cost of that placement. In this paper, we show that this
myopic strategy satisfies an analogous asymptotic optimality result
when the convergence of the posterior distribution is considered as
a function of the total cost (as opposed to the number of
observations).
\end{abstract}

%
\end{frontmatter}


\section{Introduction}

The theoretical framework of this paper is that of Bayesian adaptive
estimation with an information based objective function (see, e.g.,
MacKay \cite{mackay1992ch4}, Kujala and Lukka \cite{kujalalukka2006},
Kujala \cite{kujala2011review}). Following
the notation of Kujala \cite{bestvalue,kujala2011review}, the basic problem
we consider is the estimation of an unobservable random variable
$\Theta\dvtx \Omega\mapsto\mathsfTheta $ based on a sequence
$y_{x_1},\ldots,y_{x_t}$ of independent (given~$\theta$) realizations
from some conditional densities $p(y_{x_t}\mid\theta)$ indexed by
trial \emph{placements} $x_t$, each of which can be adaptively chosen
from some set $\mathsf{X}$ based on the outcomes
$(y_{x_1},\ldots,y_{x_{t-1}})$ of the earlier observations. A commonly
used greedy strategy is to choose the next placement so as to maximize
the expected immediate information gain, that is, the decrease of the
(differential) entropy of the posterior distribution given the next
observation.

Previous work on the asymptotics of Bayesian estimation (see, e.g.,
Schervish \cite{schervish1995}, van~der Vaart \cite{vandervaart1998}) has mostly concentrated on
the i.i.d. case, and in the few cases where the independent (given
$\theta$) but not identical case is considered, it is customarily
assumed that a certain fixed sequence of variables is given. Hence,
these results do not apply to the present situation where the sequence
$X_t$ of placements is also random.

Paninski \cite{paninski2005} has developed an asymptotic theory for this
adaptive setting. He states consistency and asymptotic normality
results for the greedy information maximization placement strategy and
quantifies the asymptotic efficiency of the method. However, the
proofs therein are not complete and hence do not provide a sufficient
foundation for some generalizations and theorems we are interested in.
In this paper, we develop a more general theory which allows us to
generalize the main results of Paninski \cite{paninski2005} to almost sure
convergence (with novel proofs) and to show that the greedy method is
in a certain sense asymptotically optimal among \emph{all} placement
methods. Furthermore, we provide a rigorous and general framework
that lends itself to further extensions of the theory.

One particular extension we are interested in is analyzing the
asymptotic properties of the novel framework proposed in
Kujala \cite{bestvalue}. In this framework, the observation of $Y_x$ is
associated with some random cost $C_x$ (see Section~\ref{secvaryingcost}
for details). To make measurement ``cost-effective'', a myopic
placement rule is considered that on each trial $t$ maximizes the
expected value of the information gain (decrease of entropy)
\[
G_t=\mathrm{H}(\Theta\mid Y_{X_1},\ldots, Y_{X_{t-1}})-
\mathrm {H}(\Theta\mid Y_{X_1},\ldots, Y_{X_t})
\]
divided by the expected value of the cost $C_t=C_{X_t}$. This is
called a myopic strategy as it looks only one step ahead. However, it
is not a greedy strategy as it does not optimize the immediate gain.

In Kujala \cite{bestvalue}, the following fairly simple asymptotic optimality
result is given for this myopic strategy.

\begin{theorem}\label{asymptotic2}
Suppose that there exists a constant $\alpha>0$ such that
%
%
\begin{equation}
\label{const} \max_{x\in\mathsf{X}}\frac{\mathrm{E}(G_t\mid\mathbf{y}, X_t=x)}{\mathrm{E}(C_t\mid\mathbf{y}, X_t=x)} = \alpha
\end{equation}
for all possible sets $\mathbf{y}$ of past observations. If the
next placement $X_t$ is defined as the maximizer of (\ref{const})
and if for some $\sigma^2<\infty$ and $\varepsilon>0$,
%
%
\begin{equation}
\label{ass1} \cases{\displaystyle \Var(G_t\mid Y_{X_1},\ldots,Y_{X_{t-1}}) \le\sigma^2,
\vspace*{3pt}\cr
\displaystyle
\Var(C_t\mid Y_{X_1},\ldots,Y_{X_{t-1}}) \le
\sigma^2,
\vspace*{3pt}\cr
\displaystyle \mathrm{E}(C_t\mid
Y_{X_1},\ldots,Y_{X_{t-1}}) \ge \varepsilon}
\end{equation}
for all $t$, then the gain-to-cost ratio satisfies
\[
\lim_{t\to\infty} \frac{G_1+\cdots+G_t}{
C_1+\cdots+C_t}\mathop{=}\limits
^{\mathit{a.s.}} \alpha.
\]
This is asymptotically optimal in the sense that for any other
strategy that satisfies (\ref{ass1}), we have
\[
\limsup_{t\to\infty} \frac{G_1+\cdots+G_t} {C_1+\cdots+C_t} \mathop{\le}\limits
^{\mathit{a.s.}} \alpha.
\]
\end{theorem}

However, this result requires the obtainable information gains to not
decrease over time for the optimality condition to make sense and
hence does not in general apply to smooth models. In this paper, we
provide a counterpart of the above result using an optimality
criterion (D-optimality) relevant to smooth models.

Our results are structured as follows. In
Section~\ref{secconsistency}, we derive strong consistency of the
posterior distributions under extremely mild, purely topological
conditions on the family of likelihood functions. In
Section~\ref{secnormality}, we consider the local smoothness
assumptions (to be assumed in a certain neighborhood of the true
parameter value) required for asymptotic normality. In
Section~\ref{secasymprop}, we develop a theory of asymptotic
proportions and use it for a novel type of convergence of random
variables that is required in our analysis. Then, in
Sections~\ref{secd-optimality} and \ref{secasyment}, we are able to
quantify the asymptotic covariance and asymptotic entropy of the
posterior distribution and to show a form of asymptotic optimality for
the standard greedy information maximization strategy. In
Section~\ref{secvaryingcost}, these results are generalized to the
situation with random costs of observation associated with each
placement as discussed above. The heuristically justified, myopic
placement strategy proposed in Kujala \cite{bestvalue} turns out to be
asymptotically optimal also in the sense of the present paper,
supporting the view that this strategy is the most natural
generalization of the greedy information maximization strategy to the
situation where the costs of observation can vary. We give concrete
examples of the optimality results in Section~\ref{secexamples} and
then end with general discussion in Section~\ref{secdiscussion}.

\subsection{Preliminaries}

We shall denote random variables by upper case letters and their
specific values by lower case letters. The information theoretic
definitions that we will use are the (differential) entropy $\mathrm
{H}(A) =
-\int p(a)\log p(a) \,\mathrm{d}a$, which does depend on the parameterization of
$a$, the Kullback--Leibler divergence
\[
D_{\KL} \bigl(p(a) \,\Vert\, p(b) \bigr) = \int p(a)\log
\frac{p(a)}{p(b)}\,\mathrm{d}a,
\]
which is independent of the parameterization, and the mutual
information
\begin{eqnarray*}
\I(A;B) &=& \int p(a,b)\log\frac{p(a,b)}{p(a)p(b)}\,\mathrm{d}(a,b)
\\
&=& \int p(a) D_{\KL} \bigl(p(b\mid a) \,\Vert\, p(b) \bigr)\,\mathrm{d}a
\\
&=& \int p(b) D_{\KL} \bigl(p(a\mid b) \,\Vert\, p(a) \bigr)\,\mathrm{d}b,
\end{eqnarray*}
which is also independent of the parameterization as well as
symmetric. Also, the identities $\I(A;B) = \mathrm{H}(A) - \mathrm
{E}(\mathrm{H}(A\mid
B)) = \mathrm{H}(B) - \mathrm{E}(\mathrm{H}(B\mid A))$ hold whenever
the differences are
well defined. This is all standard notation (see, e.g.,
Cover and Thomas \cite{coverthomas2006}) except that in our notation, there is no
implicit expectation over the values of $A$ in $\mathrm{H}(B\mid A)$,
and so
it is a random variable depending on the value of $A$. Similarly, a
conditional density $p(b\mid a)$ as an argument to $D_{\KL}(\cdots)$ is
treated the same way as any other density of $b$, with no implicit
expectation over $a$.

The densities $p(a)$ and $p(b)$ above are assumed to be taken
w.r.t. arbitrary dominating measures ``$\mathrm{d}a$'' and ``$\mathrm{d}b$''. Thus,
following Lindley \cite{lindley1956}, we are in fact working in full measure
theoretic generality even though we use the more familiar notation.
The underlying probability space is $(\Omega,\mathcal{F},\Pr)$ and so,
for example, $\Pr\{\Theta\in U\}$ means the probability that the
value of
$\Theta\dvtx \Omega\to\mathsfTheta $ is within the measurable set
$U\subset\mathsfTheta $. In some places we may abbreviate this by
$p(U)$, but it will be clear from the context what random variable is
referred to. When we say ``for a.e. $\theta$'', it is w.r.t. the
prior distribution of $\Theta$. The $\sigma$-algebra of
$\mathsfTheta $ is assumed to contain at least the Borel sets of the
topology which $\mathsfTheta $ is assumed to be endowed with.

For any fixed $x\in\mathsf{X}$, we assume that the conditional
densities $p(y_x\mid\theta)$ are given w.r.t. the same dominating
$\sigma$-finite measure ``$\mathrm{d}y_x$'' for all $\theta\in\mathsfTheta $
and when we say ``for a.e. $y_x$'', it is w.r.t. this measure.
For brevity, we shall indicate conditioning on the data $\mathbf{Y}_t:= (Y_{X_1},\ldots,Y_{X_t})$ by the subscript $t$ on any quantities
that depend on them. For example, $p_t(\theta) = p(\theta\mid
\mathbf{Y}_t)$ is the posterior density of $\Theta$ given
$\mathbf{Y}_t$ and $\mathrm{E}_t(f(\Theta)) = \mathrm{E}_t(f(\Theta
)\mid
\mathbf{Y}_t)$ is the posterior expectation of $f(\Theta)$ given~$\mathbf{Y}_t$.

It is often assumed that one can observe multiple independent (given
$\theta$) copies of the same random variable $Y_x$. However, instead
of complicating the general notation with something like
$Y_{x_t}^{(t)}$, we rely on the fact that the set $\mathsf{X}$ can
explicitly include separate indices for any identically distributed
copies, for example, one might have $[Y_{(x,t)}\mid\theta]
\mathop{\sim}\limits^{\mathrm{i.i.d.}}\, [Y_{(x,t')}\mid\theta]$
for all
$t,t'\in\mathbb{N}$, $t\ne t'$. Hence, we can use the simple notation
with no loss of generality.

The greedy information gain maximization strategy can be formally
defined as choosing the placement $X_t$ to be the value $x$ that
maximizes the mutual information $\I_{t-1}(\Theta;Y_x) =
\mathrm{H}_{t-1}(\Theta) - \mathrm{E}_{t-1}(\mathrm{H}_{t-1}(\Theta
\mid Y_x))$, the
expected decrease in the entropy of $\Theta$ after the next
observation. In some models, there may be no maximum of the mutual
information in which case the placement should be chosen sufficiently
close to the supremum, which we formally define as the ratio of the
mutual information and its supremum converging to one (condition O4 in
Section~\ref{secoptimality}).

%
\section{Consistency}\label{secconsistency}

The general assumptions for consistency are:
\begin{enumerate}[C3.]
\item[C1.] The parameter space $\mathsfTheta $ is a compact
topological space.
\item[C2.] The family of log-likelihoods is (essentially)
equicontinuous, that is, for all $\theta\in\mathsfTheta $ and
$\varepsilon>0$, there exists a neighborhood $U$ of $\theta$ such that
whenever $\theta'\in U$,
\[
\bigl\llvert \log p(y_x\mid\theta)-\log p \bigl(y_x
\mid \theta' \bigr) \bigr\rrvert < \varepsilon
\]
for a.e. $y_x$ for all $x\in\mathsf{X}$.
\item[C3.] All points in $\mathsfTheta $ are statistically
distinguishable from each other. That is, for all distinct
$\theta,\theta'\in\mathsfTheta $,
\[
d_x \bigl(\theta,\theta' \bigr):= \int \bigl\llvert
p(y_x\mid\theta)-p_t \bigl(y_x\mid\theta
' \bigr) \bigr\rrvert \,\mathrm{d}y_x > 0
\]
for some $x\in\mathsf{X}$.
\item[C4.] For some $\gamma>0$, the placements $X_t$ satisfy
\[
\I_{t-1}(\Theta;Y_{X_t})\ge\gamma\sup_{x\in\mathsf{X}}
\I _{t-1}(\Theta;Y_x)
\]
for all sufficiently large $t$.
\end{enumerate}
%

\begin{remark}
These assumptions for consistency are considerably weaker than those
formulated in Paninski \cite{paninski2005}. In particular, the assumptions
\textup{C1--C3} only pertain to the likelihood function $p(y_x\mid\theta)$,
absolutely nothing is assumed about the prior distribution of
$\Theta$. Furthermore, these assumptions are purely topological in
the sense that they are preserved by all homeomorphic transformations
of $\mathsfTheta $. Also, in C4, we do not require perfect
maximization of information gain; this is useful as it allows us to
apply the same result to the non-greedy strategy discussed in
Section~\ref{secvaryingcost} as well.
\end{remark}


\begin{remark}
Non-compact spaces can be handled if the log-likelihood has an
(essentially) equicontinuous extension to a compactification of
$\mathsfTheta $. This happens precisely when the following
conditions hold:
\begin{enumerate}[C3$'$.]
\item[C1$'$.] The parameter space $\mathsfTheta $ is a topological
space.
\item[C2$'$.] The function $f(\theta)=((x,y_x)\mapsto\log
p(y_x\mid\theta))$, with the topology of the target space induced
by the ($[0,\infty]$-valued) norm
\[
\llVert v\rrVert = \sup_{x\in\mathsf{X}}\, \mathop{\operatorname{ess}\operatorname{sup}}_{y_x} \bigl\llvert v(x,y_x) \bigr\rrvert,
\]
is continuous (this is just restating C2) and the closure of the
range $f(\mathsfTheta )$ is compact (this is the extra condition
needed for non-compact spaces).
\item[C3$'$.] For all distinct $\theta,\theta'\in\mathsfTheta $,
the inequality $f(\theta) \ne f(\theta')$ holds true, where
equality is interpreted w.r.t. a.e. $y_x$. (This is equivalent
to C3.)
\end{enumerate}
In that case, $f$ lifts continuously to the Stone--{\v C}ech
compactification $\beta\mathsfTheta $ of $\mathsfTheta $
(Theorem~\ref{thmstonecech}). Condition~\textup{C3} may not hold for the points added
by the compactification, but this can be fixed by moving to the
compact quotient space $\beta\mathsfTheta /\ker(f)$. Thus, \textup{C1--C3}
can always be replaced by the strictly weaker conditions
\textup{C1$'$--C3$'$}.
\end{remark}

%
\begin{lemma}\label{lemmingain}
Suppose that \textup{C1--C3} hold. Then, there exists a metric
$d\dvtx \mathsfTheta \times\mathsfTheta \to\mathbb{R}$ that is
consistent with the topology of $\mathsfTheta $, and an estimator
$\hat\Theta_t$ such that for each $t$ there exists $x\in\mathsf{X}$
such that
\[
\I_t(Y_x;\Theta) \ge\mathrm{E}_t \bigl(d(
\Theta,\hat\Theta_t)^2 \bigr).
\]
\end{lemma}

\begin{pf}
First, we show that the pseudometric $d_x$ defined in C3 is
continuous in $\mathsfTheta \times\mathsfTheta $ for all
$x\in\mathsf{X}$. It can be shown using C2 that for any
$\theta\in\mathsfTheta $ and $\varepsilon>0$, there exists a neighborhood
$U_{\theta,\varepsilon}$ such that $d_x(\theta,\theta')\le
\varepsilon$ for all
$\theta'\in U_{\theta,\varepsilon}$. Thus, for any $\varepsilon>0$ and
$\theta_1,\theta_2\in\mathsfTheta $, the triangle inequality
implies
\[
\bigl\llvert d_x \bigl(\theta_1',
\theta'_2 \bigr) - d_x(
\theta_1,\theta_2) \bigr\rrvert \le d_x
\bigl( \theta_1,\theta_1' \bigr) +
d_x \bigl(\theta_2,\theta_2'
\bigr) \le2 \varepsilon
\]
whenever $(\theta_1',\theta_2') \in U_{\theta_1,\varepsilon}\times
U_{\theta_2,\varepsilon}$, and so $d_x$ is continuous.

As $d_x$ is continuous, the set
\[
S_x = \bigl\{ \bigl(\theta,\theta' \bigr)\in\mathsfTheta \times\mathsfTheta\dvt d_x \bigl(\theta,\theta'
\bigr)>0 \bigr\}
\]
is open for every $x\in\mathsf{X}$. Now C3 implies that
$\bigcup_{x\in\mathsf{X}}S_x$ covers
$\mathsfTheta \times\mathsfTheta $, and as
$\mathsfTheta \times\mathsfTheta $ is compact, there exists a
finite subcover $\bigcup_{x\in\mathsf{X}'}S_x$. It follows that
\[
d \bigl(\theta,\theta' \bigr) = \biggl[\frac{1}{8\llvert  \mathsf{X}'\rrvert  }\sum
_{x\in
\mathsf{X}'} \biggl(\int \bigl\llvert p(y_x\mid
\theta)-p_t \bigl(y_x\mid\theta' \bigr)
\bigr\rrvert \,\mathrm{d}y_x \biggr)^2 \biggr]^{1/2}
\]
is positive definite and hence a metric. Since $\mathsf{X}'$ is
finite, this metric inherits the continuity of $d_x$.

To show that the topology induced by $d$ coincides with that of
$\mathsfTheta $, let $U$ be an arbitrary open neighborhood of
$\theta_0$. Then $U^c$ is compact and so its continuous image $S:=
\{ d(\theta_0,\theta)\dvt \theta\in U^c \}$ is compact, too. It
follows that $S^c$ is open and as $0\in S^c$, we obtain
$[0,\delta_U)\subset S^c$ for some $\delta_U>0$. Thus, we obtain
$\{ \theta\in\mathsfTheta \dvt d(\theta_0,\theta)<\delta_U \}
\subset
U$, and so the topology induced by $d$ is finer than the default
topology of $\mathsfTheta $. As $d$ is continuous, we obtain the
converse, and so the topologies coincide.

Let then $t$ be arbitrary. We extend $d(\theta,\theta')$ with a
special point $\bar\Theta_t\notin\mathsfTheta $ for which we
define the distances
\[
d(\theta,\bar\Theta_t) = \biggl[\frac{1}{8\llvert  \mathsf{X}'\rrvert  }\sum
_{x\in\mathsf{X}'} \biggl(\int \bigl\llvert p(y_x\mid\theta
)-p_t(y_x) \bigr\rrvert \,\mathrm{d}y_x
\biggr)^2 \biggr]^{1/2}.
\]
The extended distance function may not be strictly positive
definite, but it is still a pseudometric and satisfies the triangle
inequality. Denoting
\[
\hat\Theta_t = \argmin_{\theta\in\mathsfTheta } d(\theta,\bar
\Theta_t),
\]
we\vspace*{1pt} have $d(\theta,\bar\Theta_t)\ge d(\hat\Theta_t,\bar\Theta_t)$ for
all $\theta\in\mathsfTheta $, and the triangle inequality yields
$d(\theta,\bar\Theta_t) \ge d(\theta,\hat\Theta_t) -
d(\hat\Theta_t,\bar\Theta_t)$. Adding both inequalities, we obtain
$2d(\theta,\bar\Theta_t) \ge d(\theta,\hat\Theta_t)$ for all
$\theta\in\mathsfTheta $. Now, the $L^1$-bound of
Kullback--Leibler divergence \cite{coverthomas2006}, Lemma~11.6.1,
yields
%
%
\begin{eqnarray*}
\max_{x\in\mathsf{X}'} \I_t(Y_x;\Theta) &\ge&
\frac{1}{\llvert  \mathsf{X}'\rrvert  }\sum_{x'\in\mathsf{X}'} \I_t(Y_x;
\Theta)
\\
&=& \int\frac{1}{\llvert  \mathsf{X}'\rrvert  }\sum_{x'\in\mathsf{X}'} D_{\KL}
\bigl( p(y_x\mid\theta) \,\Vert\, p_t(y_x)
\bigr)p_t(\theta)\,\mathrm{d}\theta
\\
{ \bigl(L^1\mbox{ bound} \bigr)}\qquad &\ge&\int\frac{1}{\llvert  \mathsf{X}'\rrvert  }\sum
_{x'\in\mathsf{X}'}\frac{1}2 \biggl[\int \bigl\llvert
p(y_x \mid\theta)-p_t(y_x) \bigr\rrvert
\,\mathrm{d}y_x \biggr]^2p_t(\theta )\,\mathrm{d}\theta
\\
&=& 4\int d(\theta,\bar\Theta_t)^2p_t(
\theta)\,\mathrm{d}\theta \ge\int d(\theta,\hat\Theta_t)^2p_t(
\theta)\,\mathrm{d}\theta.
\end{eqnarray*}\upqed
\end{pf}


\begin{lemma}\label{lemfinitesum}
Suppose that $K$ is a function of $\Theta$ and has a finite range
$\mathsf{K}$. Then, for arbitrarily chosen placements $X_t$, the
inequality $\sum_{t=1}^\infty\I_{t-1}(K; Y_{X_t}) < \infty$ holds
almost surely (which implies $\I_{t-1}(K; Y_{X_t})\astoi 0$).
\end{lemma}

\begin{pf}
As $\I_{t-1}(K;Y_{X_t}) =
\mathrm{H}_{t-1}(K)-\mathrm{E}_{t-1}(\mathrm{H}_t(K))$, where $0\le
\mathrm{H}_t(K)\le\log\llvert  \mathsf K\rrvert  $ for all $t$, we obtain
\begin{eqnarray*}
\mathrm{E} \Biggl( \sum_{k=1}^t
\I_{k-1}(K;Y_{X_k}) \Biggr) &=& \mathrm{E} \bigl(
\mathrm{H}_0(K) - \mathrm{E}_{t-1} \bigl(\mathrm
{H}_t(K) \bigr) \bigr) \le\log\llvert \mathsf K\rrvert
\end{eqnarray*}
for all $t$. As $\I_{t-1}(K;Y_{X_t})$ is nonnegative, the
sequence of partial sums is non-decreasing, and Lebesgue's monotone
convergence theorem yields
\[
\mathrm{E} \Biggl( \sum_{k=1}^\infty
\I_{k-1}(K;Y_{X_k}) \Biggr) = \lim_{t\to\infty}
\mathrm{E} \Biggl( \sum_{k=1}^t \I
_{k-1}(K;Y_{X_k}) \Biggr)\le\log\llvert \mathsf K\rrvert <
\infty,
\]
which implies the statement.
\end{pf}


\begin{lemma}\label{leminfoconv}
Suppose that \textup{C1} and \textup{C2} hold. Then $\I_{t-1}(\Theta; Y_{X_t})\astoi0$
for arbitrarily chosen placements $X_t$.
\end{lemma}

\begin{pf}
Let $\varepsilon> 0$ be arbitrary. As $\mathsfTheta $ is compact, a
finite number of the sets $U_{\theta,\varepsilon}$ given by C2 cover it.
Thus, we can partition the parameter space into a finite number of
subsets $\mathsfTheta _k$ each one contained in some
$U_{\theta,\varepsilon}$. Letting the random variable $K$ denote the index
of the subset that $\Theta$ falls into, the chain rule of mutual
information yields
%
%
\begin{equation}
\label{infochain} \I_{t-1}(\Theta;Y_t) = \I_{t-1}(
\Theta,K;Y_t) = \I_{t-1}(K;Y_t) + \sum
_k p_{t-1}(k)\I_{t-1}(
\Theta;Y_t\mid k),
\end{equation}
where $Y_t:= Y_{X_t}$ and Lemma~\ref{lemfinitesum} implies that
$\I_{t-1}(K;Y_t)\asto0$. Let us then look at the latter term.
Convexity of the Kullback--Leibler divergence yields
\begin{eqnarray*}
\I_{t-1}(\Theta;Y_t\mid k) &=& \int p_{t-1}(\theta
\mid k)D_{\KL} \bigl(p(y_t\mid\theta) \,\Vert\,
p_{t-1}(y_t\mid k) \bigr)\,\mathrm{d}\theta
\\
&\le&\int p_{t-1}(\theta\mid k) \biggl[\int p_{t-1} \bigl(
\theta'\mid k \bigr)D_{\KL} \bigl(p(y_t\mid
\theta) \,\Vert\, p \bigl(y_t\mid\theta' \bigr) \bigr)\,\mathrm{d}
\theta ' \biggr]\,\mathrm{d}\theta
\\
&=& \int\!\!\!\int p_{t-1}(\theta\mid k)p_{t-1} \bigl(
\theta'\mid k \bigr) \biggl[\int p(y_t\mid\theta)
\underbrace{\log\frac{p(y_t\mid\theta
)}{p(y_t\mid\theta')}}_{\le2\varepsilon
\mathrm{~for~a.e.~}y_t}\,\mathrm{d}y_t \biggr]\,\mathrm{d}\theta \,\mathrm{d}
\theta'
\\
&\le&2\varepsilon
\end{eqnarray*}
for all $t$. Thus,
\[
\limsup_{t\to\infty}\I_{t-1}(\Theta;Y_t)\le2
\varepsilon
\]
almost surely. As $\varepsilon>0$ was arbitrary, we obtain
$\I_{t-1}(\Theta;Y_t)\asto0$.
\end{pf}


\begin{lemma}\label{lemeconv}
For any measurable function $f\dvtx \mathsfTheta \to\mathbb{R}$, if the
prior expectation $\mathrm{E}f(\Theta)$ is well-defined and finite, then
$\lim_{t\to\infty}\mathrm{E}_t f(\Theta)$ exists as a finite
number almost
surely.
\end{lemma}

\begin{pf}
The finiteness of $\mathrm{E}f(\Theta)$ implies that $\mathrm
{E}\llvert  f(\Theta)\rrvert  $ must
also be finite and so $Z_t:= \mathrm{E}_t f(\Theta)$ satisfies
$\mathrm{E}\llvert  Z_t\rrvert   =
\mathrm{E}\llvert  \mathrm{E}_tf(\Theta)\rrvert   \le\mathrm{E}\llvert  f(\Theta)\rrvert  <\infty
$ for all $t$.
Furthermore, since $Z_{t+1}$ depends linearly on the posterior
$p_{t+1}$ whose expectation $\mathrm{E}_t(p_{t+1})$ equals the prior $p_t$,
we obtain $\mathrm{E}_t(Z_{t+1}) = Z_t$ for all $t$ and so $Z_t$ is a
martingale. As $\sup_t\mathrm{E}\llvert  Z_t\rrvert   \le\mathrm{E}\llvert  f(\Theta
)\rrvert  <\infty$,
Theorem~\ref{thmmartingaleconv} implies that $\lim Z_t$ exists as a
finite number almost surely.
\end{pf}


\begin{theorem}[(Strong consistency)]\label{thmsc}
Suppose that \textup{C1--C4} hold. Then, conditioned on almost any
$\theta_0\in\mathsfTheta $ as the true parameter value, the
posteriors are strongly consistent, that is, $\Pr_t\{\Theta\in
U\}\astoi1$ for any neighborhood $U$ of $\theta_0$.
\end{theorem}

\begin{pf}
As the metric $d$ given by Lemma~\ref{lemmingain} is bounded,
Lemma~\ref{lemeconv} implies that
$\lim_{t\to\infty}\mathrm{E}_t(d(\Theta,\theta))$ exists and is
finite for
all $\theta$ in a countable dense subset of $\mathsfTheta $ almost
surely, in which case continuity of $d$ implies the same for all
$\theta\in\mathsfTheta $.

Lemmas~\ref{lemmingain} and \ref{leminfoconv} and C4 yield
$\mathrm{E}_t(d(\Theta,\hat\Theta_t))\asto0$. As $d$ is bounded, Lebesgue's
dominated convergence theorem and Markov's inequality imply
\[
\Pr \bigl\{d(\Theta,\hat\Theta_t)>\varepsilon \bigr\} \le
\frac{\mathrm{E}(d(\Theta,\hat\Theta_t))}{\varepsilon} = \frac{\mathrm{E}(\mathrm{E}_t(d(\Theta,\hat\Theta
_t)))}{\varepsilon} \to0
\]
for all $\varepsilon>0$ and so $d(\Theta,\hat\Theta_t)\pto0$. Convergence
in probability implies that there exists a subsequence $t_k$ such
that $d(\Theta,\hat\Theta_{t_k})\asto0$. Thus, conditioned on
almost any $\theta_0$ as the true value, we obtain
$d(\theta_0,\hat\Theta_{t_k})\asto0$, and the triangle inequality
yields
\[
\mathrm{E}_{t_k} \bigl(d(\Theta,\theta_0) \bigr) \le
\mathrm{E}_{t_k} \bigl(d(\Theta,\hat\Theta_{t_k}) \bigr) +d(
\theta_0,\hat\Theta_{t_k})\asto0.
\]
As we have already established that the full sequence
$\mathrm{E}_t(d(\Theta,\theta_0))$ almost surely converges, it now follows
that the limit must almost surely be zero. Thus, given any
neighborhood $U\supset B_d(\theta_0,\varepsilon)$ of $\theta_0$, Markov's
inequality yields
\[
\Pr_t \bigl\{\Theta\in U^c \bigr\}\le\Pr_t
\bigl\{\Theta\in B_d(\theta _0,\varepsilon)^c
\bigr\}\le \frac{\mathrm{E}_t(d(\Theta,\theta_0))}{\varepsilon}\asto0.
\]\upqed
%
%
\end{pf}



\begin{lemma}\label{lemscimpl}
Suppose that \textup{C1--C3} hold and assume that conditioned on
$\theta_0\in\mathsfTheta $ as the true parameter value, the
posteriors are strongly consistent. Then:
\begin{enumerate}
\item Given any metric $d$ consistent with the topology of
$\mathsfTheta $,
\[
\Theta^*_t:= \argmin_{\theta\in\mathsfTheta }\mathrm {E}_t
\bigl(d(\Theta,\theta)^2 \bigr)\astoi\theta_0.
\]
\item For any neighborhood $U$ of $\theta_0$ there exists a constant
$c>0$ such that, almost surely, $\I_t(Y_x;\Theta)\ge
c\Pr_t\{\Theta\in U^c\}$ for some $x\in\mathsf{X}$ for all
sufficiently large $t$.
\end{enumerate}
\end{lemma}

\begin{pf}
Let $D$ be the diameter of $\Theta$. The triangle inequality $a\le
b+c$ implies $a^2 \le(b+c)^2 \le2(b^2+c^2)$ and so consistency of
the posteriors yields
\begin{eqnarray*}
d \bigl(\theta_0,\Theta^*_t \bigr)^2 &\le& 2
\mathrm{E}_t \bigl(d(\Theta,\theta_0)^2+d
\bigl(\Theta,\Theta^*_t \bigr)^2 \bigr) \le 4
\mathrm{E}_t \bigl(d(\Theta,\theta_0)^2
\bigr)
\\
&\le& 4 \bigl(r^2+D^2\Pr_t \bigl\{\Theta\in
B_d(\theta_0,r)^c \bigr\} \bigr) \asto 4
\bigl(r^2+D^2\cdot0 \bigr)
\end{eqnarray*}
for all $r>0$, which implies $\Theta^*_t\asto\theta_0$.

Let us then assume that the metric $d$ is the one given by
Lemma~\ref{lemmingain} and choose $\varepsilon>0$ such that
$B_d(\theta_0,2\varepsilon)\subset U$. As $\Theta^*_t\asto\theta
_0$, we
have $B_d(\Theta^*_t,\varepsilon)\subset U$ for all sufficiently
large $t$,
and so Lemma~\ref{lemmingain} and Markov's inequality yield
\begin{eqnarray*}
\I_t(Y_x;\Theta)&\ge&\mathrm{E}_t \bigl(d(
\Theta,\hat\Theta_t)^2 \bigr)
\\
&\ge& \mathrm{E}_t \bigl(d \bigl(\Theta,\Theta^*_t
\bigr)^2 \bigr) \ge\varepsilon^2 \Pr_t \bigl\{
\Theta\in B_d \bigl(\Theta^*_t,\varepsilon
\bigr)^c \bigr\} \ge\varepsilon^2 \Pr_t \bigl
\{ \Theta\in U^c \bigr\}
\end{eqnarray*}
for some $x\in\mathsf{X}$.
\end{pf}

\subsection{Asymptotic entropy}

The differential entropy is sensitive to the parameterization, but
asymptotically, we can in most cases ignore this due to the following
lemma.

\begin{lemma}\label{lementkleq}
Suppose that the prior entropy $\mathrm{H}(\Theta)$ is well-defined and
finite. Then,
\[
\lim_{t\to\infty} \bigl[ \mathrm{H}_t(\Theta) +
D_{\KL} \bigl(p_t(\theta ) \,\Vert\, p(\theta) \bigr) \bigr]
\]
exists as a finite number almost surely.
\end{lemma}

\begin{pf}
As $\mathrm{H}_t(\Theta) + D_{\KL}(p_t(\theta) \,\Vert\, p(\theta)) =
\mathrm{E}_t
\log p(\Theta)$ and $\mathrm{E}\log p(\Theta) = -\mathrm{H}(\Theta
)$ is
well-defined and finite, the statement follows from
Lemma~\ref{lemeconv}.
\end{pf}


\begin{lemma}\label{lemmaxchange}
Suppose that \textup{C1$'$} holds and let $f$ be defined as in \textup{C2$'$}. Then,
for any subset $S\subset\mathsfTheta $,
\[
\bigl\llvert \log p_{t+1}(\theta\mid S)-\log p_t(\theta\mid
S) \bigr\rrvert \le 2 \diam f(S)
\]
for all $\theta\in S$. If \textup{C2$'$} holds, then this upper bound is
finite.
\end{lemma}

\begin{pf}
Let $\theta_1\in S$ be fixed. If $p_t(\theta\mid S)$ is multiplied
by $p(y_x\mid\theta)/p(y_x\mid\theta_1)$, it can change by at most a
factor of $\exp(\diam f(S))$, and for the same reason, the
normalization constant for this density is within a factor of
$\exp(\diam f(S))$ from $1$. The statement follows.

Suppose then that C2$'$ holds. As $f(\mathsfTheta )$ is compact,
it follows that $f(S)\subset f(\mathsfTheta )$ must be bounded.
\end{pf}


\begin{lemma}\label{lementsublin}
Suppose that \textup{C1} and \textup{C2} hold. Then, for any $\varepsilon>0$, the inequality
$D_{\KL}(p_t(\theta) \,\Vert p(\theta))<\varepsilon t$ holds true for all
sufficiently large $t$.
\end{lemma}

\begin{pf}
Let $\varepsilon>0$ be arbitrary. As in the proof of
Lemma~\ref{leminfoconv}, we partition $\mathsfTheta $ into a
finite number of subsets $\mathsfTheta _k$ such that $\llvert  \log
p(y_x\mid\theta)-\log p(y_x\mid\theta_k)\rrvert  \le\varepsilon$ for all
$\theta\in\mathsfTheta _k$, $y_x$, and $x\in\mathsf{X}$, where
$\theta_k$ is some fixed point of $\mathsfTheta _k$. Let the
random variable $K$ denote the index of the subset that $\Theta$
falls into. Lemma~\ref{lemmaxchange} implies that
\[
\bigl\llvert \log p_{t+1}(\theta\mid k) - \log p_t(\theta
\mid k) \bigr\rrvert \le2 \varepsilon
\]
for all $\theta\in\mathsfTheta _k$, which yields
\[
D_{\KL} \bigl(p_t(\theta\mid k) \,\Vert\, p(\theta\mid k)
\bigr) = \mathrm{E}_t \biggl(\log\frac{p_t(\Theta\mid k)}{p(\Theta\mid k)} \Bigm| k \biggr)
\le2 \varepsilon t
\]
for all $t$ and $k$. The chain rule of Kullback--Leibler divergence
now yields
\begin{eqnarray*}
D_{\KL} \bigl(p_t(\theta) \,\Vert\, p(\theta) \bigr) &=&
D_{\KL} \bigl(p_t(k) \,\Vert\, p(k) \bigr) + \sum
_k p_t(k)D_{\KL} \bigl(p_t(
\theta\mid k) \,\Vert\, p(\theta\mid k) \bigr)
\\
&\le&\log\max_k p(k)^{-1} + 2\varepsilon t,
\end{eqnarray*}
where we may assume that $p(k)$ is positive since we can drop any
set $\mathsfTheta _k$ with $p(k)=0$ from the partition.
\end{pf}

\iftrue

\begin{lemma}\label{leminfoupper}
Suppose that $\mathsfTheta \subset\mathbb{R}^n$ is bounded and the
family of log-likelihoods is uniformly Lipschitz, that is,
\[
\bigl\llvert \log p(y_x\mid\theta)-\log p \bigl(y_x
\mid \theta' \bigr) \bigr\rrvert \le M \bigl\llvert \theta-
\theta' \bigr\rrvert
\]
for all $\theta,\theta'\in\mathsfTheta $ for all $y_x$ and
$x\in\mathsf{X}$. Then, for arbitrarily chosen placements $X_t$,
the expected gain over $t$ trials is bounded by
$\I(\Theta;\mathbf{Y}_t)\le n\log t + c$ for some constant
$c<\infty$.
\end{lemma}

\begin{pf}
For each $t$, we can subdivide the bounded parameter space
$\mathsfTheta $ into $\le ct^n$ subsets $\mathsfTheta _k$, each
having diameter $\le t^{-1}$. Letting the random variable $K_t$
denote the index of the subset that $\Theta$ falls into, the chain
rule of mutual information yields
%
%
\begin{equation}
\label{eqinfoupper} \I(\Theta;\mathbf{Y}_t) = \underbrace{
\I(K_t;\mathbf{Y}_t)}_{\le
\log(ct^n)} + \sum
_{k_t} p(k_t)\underbrace{\I(\Theta;
\mathbf{Y}_t\mid k_t)}_{\le M} \le n \log t + \log
c + M
\end{equation}
as in equation~(\ref{infochain}) in the proof Lemma~\ref{leminfoconv}.
%
\end{pf}
\fi

\section{Asymptotic normality}\label{secnormality}

In this section, we assume that:
\begin{enumerate}[N3.]
\item[N1.] The parameter space $\mathsfTheta $ is 
a subset of $\mathbb{R}^n$.
\item[N2.] The true parameter value $\theta_0$ is an interior point
of $\mathsfTheta $.
\item[N3.] The\vspace*{1pt} log-likelihood $\theta\mapsto\log p(y_x\mid
\theta)$
is twice continuously differentiable with $\llvert  \nabla_\theta \log
p(y_x\mid\theta)\rrvert  \le M$ and $\llvert  \nabla^2_\theta \log
p(y_x\mid\theta)\rrvert  \le M$ for all $x\in\mathsf{X}$ and $y_x$.
\item[N4.] The family of Hessians $\theta\mapsto\nabla^2_\theta
\log
p(y_x\mid\theta)$ is equicontinuous at $\theta_0$ over all
$x\in\mathsf{X}$ and~$y_x$.
\item[N5.] The prior density is absolutely continuous w.r.t. the
Lebesgue measure with positive and continuous density at $\theta_0$.
\end{enumerate}
For simplicity of notation, all statements are implicitly conditioned
on $\theta_0$ being the true parameter value. Throughout this
section, we will denote the posterior mean and covariance by
$\hat\Theta_t:= \mathrm{E}_t(\Theta)$ and $\Sigma_t=\Cov
_t(\Theta)$. Note
that the expected square error $\mathrm{E}_t(\llvert  \Theta-\theta\rrvert  ^2)$ is minimized
by the mean $\theta= \mathrm{E}_t(\Theta)$. Thus, if the posteriors are
strongly consistent, then Lemma~\ref{lemscimpl} implies that
$\hat\Theta_t\asto\theta_0$. Note also that the square error is
related to the variance through the identity
$\mathrm{E}_t(\llvert  \Theta-\hat\Theta_t\rrvert  ^2) = \tr(\Sigma_t)$.


\begin{lemma}\label{lemiapprox}
Suppose that \textup{N1} and \textup{N3} hold and $\mathsfTheta $ is a bounded
convex set with diameter${}\le D<\infty$. Then, there exists a
constant $C_{M,D}<\infty$ such that for all $t$, and $x$,
\[
\bigl\llvert \I_t(Y_x;\Theta) - \bigl(
\tfrac{1}2\Sigma_t \bigr)\odot I_x(\hat
\Theta_t) \bigr\rrvert \le C_{M,D}\mathrm{E}_t
\bigl(\llvert \Theta-\hat\Theta_t\rrvert ^3 \bigr),
\]
where $\odot$ denotes the Frobenius product $A\odot B = \sum_{i,j}
A_{ij}B_{ij} = \tr(A^TB)$, and $I_x(\theta)$ is the Fisher
information matrix
\[
I_x(\theta):= \int \biggl[\frac{\nabla_\theta p(y_x\mid\theta)}{p(y_x\mid\theta
)} \biggr] \biggl[
\frac{\nabla_\theta p(y_x\mid\theta)}{p(y_x\mid\theta
)} \biggr]^T p(y_x\mid
\theta)\,\mathrm{d}y_x.
\]
\end{lemma}

\begin{pf}
We can formally expand the mutual information as
\begin{eqnarray*}
\I_t(Y_x;\Theta) &=& \mathrm{H}_t(Y_x)
- \mathrm{E}_t \bigl(\mathrm {H}(Y_x\mid\Theta) \bigr)
\\
&=& \int g \biggl(\int p(y_x\mid\theta)p_t(\theta)\,\mathrm{d}
\theta \biggr) \,\mathrm{d}y_x - \int \biggl(\int g \bigl(p(y_x\mid
\theta) \bigr) \,\mathrm{d}y_x \biggr)p_t(\theta )\,\mathrm{d}\theta
\\
&=& \int \biggl[ g \biggl(\int p(y_x\mid\theta)p_t(
\theta )\,\mathrm{d}\theta \biggr) - \int g \bigl(p(y_x\mid\theta) \bigr)
p_t(\theta)\,\mathrm{d}\theta \biggr] \,\mathrm{d}y_x,
\end{eqnarray*}
where $g(p) = -p\log p$. (Although $\mathrm{H}_t(Y_x) -
\mathrm{E}_t(\mathrm{H}(Y_x\mid\Theta))$ may not be well defined
here, the last
line is always well-defined and equal to the mutual information.)
Denoting $p_{y_x}:= p(y_x\mid\hat\Theta_t)$, Taylor's theorem
yields
\[
g(p) = -p_{y_x}\log p_{y_x} - (1+\log p_{y_x})
(p-p_{y_x}) - \frac
{(p-p_{y_x})^2}{2p_{y_x}} + \frac{(p-p_{y_x})^3}{6q_{p,y_x}^2},
\]
where $q_{p,y_x}$ is some number between $p_{y_x}$ and $p$. The
error term is bounded by
\[
\bigl\llvert \varepsilon_{y_x}(p) \bigr\rrvert:= \biggl\llvert
\frac
{(p-p_{y_x})^3}{6q_{p,y_x}^2} \biggr\rrvert \le\frac{\llvert  p-p_{y_x}\rrvert  ^3}{6\min\{p,p_{y_x}\}^3}p_{y_x} =
\frac{1}6 \bigl(\exp \bigl(\llvert \log p - \log p_{y_x}\rrvert
\bigr)-1 \bigr)^3 p_{y_x},
\]
and as $\llvert  \log p(y_x\mid\theta) - \log p(y_x\mid\hat\Theta_t)\rrvert  \le
M\llvert  \theta-\hat\Theta_t\rrvert  \le MD$, we further obtain
\begin{eqnarray*}
\bigl\llvert \varepsilon_{y_x} \bigl(p(y_x\mid\theta)
\bigr) \bigr\rrvert &\le& \frac{1}6 \bigl(\exp \bigl( \bigl\llvert \log
p(y_x\mid\theta) - \log p(y_x\mid\hat
\Theta_t) \bigr\rrvert \bigr)-1 \bigr)^3
p(y_x \mid\hat\Theta_t)
\\
&\le&\frac{1}6 \bigl(\exp \bigl(M\llvert \theta-\hat\Theta_t
\rrvert \bigr)-1 \bigr)^3 p(y_x\mid\hat
\Theta_t)
\\
&\le&\frac{1}6 \biggl(\frac{\exp(MD)-1}{MD}M\llvert \theta-\hat\Theta
_t\rrvert \biggr)^3 p(y_x\mid\hat
\Theta_t)
\\
&=& C_1\llvert \theta-\hat\Theta_t\rrvert ^3
p(y_x\mid\hat\Theta_t).
\end{eqnarray*}
Due to the linearity of the integral, the constant and first order
terms of the expansion cancel out, leaving just
\label{eqapprox1}
%
%
\begin{eqnarray*}
\I_t(Y_x;\Theta) &\approx& \int\frac{
- [\int p(y_x\mid\theta)p_t(\theta)\,\mathrm{d}\theta- p_{y_x} ]^2
+ \int[p(y_x\mid\theta)-p_{y_x}]^2
p_t(\theta)\,\mathrm{d}\theta
}{
2p_{y_x}
}
\,\mathrm{d}y_x
\\
&=& \int \frac{1}2\Var_t \biggl(\frac{p(y_x\mid\Theta)}{p(y_x\mid\hat\Theta
_t)} \biggr)
p(y_x\mid\hat\Theta_t)\,\mathrm{d}y_x,
\end{eqnarray*}
where the error is bounded by
%
%
\begin{eqnarray*}
&& \biggl\llvert \int\varepsilon_{y_x} \biggl(\int p(y_x
\mid \theta)p_t(\theta )\,\mathrm{d}\theta \biggr)\,\mathrm{d}_{y_x} -\int\!\!\!\int
\varepsilon_{y_x} \bigl(p(y_x\mid\theta) \bigr)
p_t(\theta)\,\mathrm{d}\theta \,\mathrm{d}y_x \biggr\rrvert
\\
&&\quad \le\int \biggl\{ \biggl\llvert \varepsilon_{y_x} \biggl(\int
p(y_x\mid \theta)p_t(\theta)\,\mathrm{d}\theta \biggr) \biggr
\rrvert +\int \bigl\llvert \varepsilon_{y_x} \bigl(p(y_x\mid
\theta) \bigr) \bigr\rrvert p_t(\theta)\,\mathrm{d}\theta \biggr\}
\,\mathrm{d}y_x
\\
&&\hspace*{-6pt}\quad \mathop{\le}\limits
^{\mathrm{Jensen}}\int \biggl\{\int \bigl\llvert \varepsilon_{y_x}
\bigl(p(y_x\mid\theta) \bigr) \bigr\rrvert p_t(\theta)\,\mathrm{d}
\theta+ \int \bigl\llvert \varepsilon_{y_x} \bigl(p(y_x\mid
\theta) \bigr) \bigr\rrvert p_t(\theta)\,\mathrm{d}\theta \biggr\}
\,\mathrm{d}y_x
\\
&&\quad \le\int2 \int C_1\llvert \theta-\hat\Theta_t\rrvert
^3 p(y_x\mid\hat\Theta _t) p_t(
\theta)\,\mathrm{d}\theta \,\mathrm{d}y_x \le2C_1 \mathrm{E}_t
\bigl(\llvert \Theta-\hat \Theta_t\rrvert ^3 \bigr)
\end{eqnarray*}
for all $t$, $\hat\Theta_t$, and $x$ (Jensen's inequality applies as
$\llvert  \varepsilon_{y_x}(p)\rrvert  $ is convex).

Now Taylor's theorem yields
\[
\frac{p(y_x\mid\theta)}{p(y_x\mid\hat\Theta_t)} = 1+ \frac{\nabla_\theta p(y_x\mid\hat\Theta_t)^T}{p(y_x\mid\hat
\Theta_t)}(\theta-\hat\Theta_t) +
\frac{1}2(\theta-\hat\Theta_t)^T
\frac{\nabla^2_\theta p(y_x\mid\theta')}{p(y_x\mid\hat\Theta_t)} (\theta-\hat\Theta_t)^T,
\]
where $\theta'$ is a convex combination of $\hat\Theta_t$ and
$\theta$. The coefficients are uniformly bounded by
\[
\biggl\llvert \frac{\nabla_\theta p(y_x\mid\hat\Theta_t)}{p(y_x\mid
\hat\Theta_t)} \biggr\rrvert = \bigl\llvert
\nabla_\theta \log p(y_x\mid\hat \Theta_t) \bigr
\rrvert \le M
\]
and
%
%
\begin{eqnarray*}
\biggl\llvert \frac{\nabla^2_\theta p(y_x\mid\theta')}{p(y_x\mid\hat
\Theta_t)} \biggr\rrvert &=& \underbrace{
\frac{p(y_x\mid\theta')}{p(y_x\mid\hat\Theta
_t)}}_{\le\exp(MD)} \bigl\llvert \underbrace{\nabla_\theta
\log p \bigl(y_x\mid\theta' \bigr)}_{\llvert  \cdot
\rrvert  \le M}
\underbrace{ \nabla_\theta \log p \bigl(y_x\mid
\theta' \bigr)^T}_{\llvert  \cdot\rrvert  \le M} + \underbrace{
\nabla^2_\theta \log p \bigl(y_x\mid
\theta' \bigr)}_{\llvert  \cdot
\rrvert  \le M} \bigr\rrvert
\\
&\le&\exp(MD) \bigl(M^2+M \bigr) =: C_2.
\end{eqnarray*}
Thus, denoting the linear term by $A$ and the error term by $B$, we
obtain
\[
\Var_t \biggl(\frac{p(y_x\mid\Theta)}{p(y_x\mid\hat\Theta
_t)} \biggr) = \Var_t(A) +
\Var_t(B) + 2\Cov_t(A,B),
\]
where
%
%
\begin{eqnarray*}
\Var_t(A) &=& \Sigma_t\odot \biggl[\frac{\nabla_\theta p(y_x\mid\hat\Theta_t)}{p(y_x\mid
\hat\Theta_t)}
\biggr] \biggl[\frac{\nabla_\theta p(y_x\mid\hat\Theta_t)}{p(y_x\mid
\hat\Theta_t)} \biggr]^T,
\\
\Var_t(B) &\le&\mathrm{E}_t \bigl(\llvert B\rrvert
^2 \bigr) \le \bigl(\tfrac{1}2 C_2
\bigr)^2 \mathrm{E}_t \bigl(\llvert \Theta-\hat
\Theta_t\rrvert ^4 \bigr) \le \bigl(\tfrac{1}2C_2
\bigr)^2 D \mathrm{E}_t \bigl(\llvert \Theta-\hat\Theta
_t\rrvert ^3 \bigr),
\\
\bigl\llvert \Cov_t(A,B) \bigr\rrvert &=& \bigl\llvert
\mathrm{E}_t(AB)-\underbrace{\mathrm {E}_t(A)}_{=0}
\mathrm{E}_t(B) \bigr\rrvert \le \mathrm{E}_t \bigl(
\llvert A\rrvert \llvert B\rrvert \bigr) \le M\tfrac{1}2C_2
\mathrm{E}_t \bigl(\llvert \Theta-\hat \Theta_t\rrvert
^3 \bigr).
\end{eqnarray*}\upqed
\end{pf}

For the next theorems and lemmas, we define the following conditions
that depend on a subset $U\subset\mathsfTheta $:
\begin{enumerate}[L3.]
\item[L1.] $\llvert  \nabla^2_\theta \log
\rrvert  p(y_x\mid\theta)-\nabla^2_\theta \log\llvert  p(y_x\mid\theta')\rrvert  <\mu/2$
for all $\theta,\theta'\in U$, $x\in\mathsf{X}$, and $y_x$.
\item[L2.] $\llvert  \log p(\theta) - \log p(\theta')\rrvert  \le C$
for all $\theta,\theta'\in U$.
\item[L3.] The maximum likelihood estimator
$\Theta^*_t:=\argmax_{\theta\in U} p(\mathbf{Y}_t\mid\theta)$ is
eventually well-defined and converges to $\theta_0$ as $t$ increases
within indices satisfying $\lambda_t\ge t\mu$, where $\lambda_t$ is
the smallest eigenvalue of $-\nabla^2_\theta \log
p(\mathbf{Y}_t\mid\theta_0)$.
\end{enumerate}

%
\begin{lemma}\label{lemL12}
Suppose that \textup{N4} and \textup{N5} hold. Then, for any $\mu,C>0$, there exists
a constant $\delta_{\mu,C}<\infty$ such that \textup{L1} and \textup{L2} hold for any
neighborhood $U$ of $\theta_0$ having diameter less than
$\delta_{\mu,C}$.
\end{lemma}

%
\begin{lemma}\label{lemoutprob}
Suppose that \textup{N1}, \textup{N3}, and \textup{L1} hold. If $p(\mathbf{Y}_t\mid\theta)\ge
p(\mathbf{Y}_t\mid\theta_0)$ for some $\theta\in U$, then
\[
\llvert \theta-\theta_0\rrvert \le\frac{2\llvert  A_t\rrvert  }{t^{1/2}\mu},
\]
where $A_t = t^{-1/2}\nabla\log p(\mathbf{Y}_t\mid\theta_0)$.
Furthermore, conditioned on $\theta_0$ as the true parameter value,
\[
\Pr \bigl\{\llvert A_t\rrvert \ge a \bigr\}\le2n\exp \biggl(-
\frac{a^2}{2nM^2} \biggr)
\]
for all $t$ satisfying $\lambda_t\ge t\mu$, where $\lambda_t$ is the
smallest eigenvalue of $-\nabla^2_\theta \log
p(\mathbf{Y}_t\mid\theta_0)$.
\end{lemma}

\begin{pf}
Taylor's theorem yields
%
%
\begin{eqnarray*}
\log p(\mathbf{Y_t}\mid\theta) &=& \log p(\mathbf{Y_t}
\mid\theta _0) + \overbrace{\nabla_\theta\log p(
\mathbf{Y_t}\mid\theta _0)}^{=:Z_t}{}^T(
\theta-\theta_0)
\\
&&{}+\underbrace{\tfrac{1}2(\theta-\theta_0)^T
\nabla_\theta^2 \log p \bigl(\mathbf{Y_t}\mid
\theta' \bigr) (\theta-\theta_0)}_{\le-(1/2)\lambda
_t\llvert  \theta-\theta_0\rrvert  ^2\le-(1/2)t\mu\llvert  \theta-\theta_0\rrvert  ^2},
\end{eqnarray*}
for some $\theta'$ between $\theta_0$ and $\theta$. Thus,
$p(\mathbf{Y}_t\mid\theta)\ge p(\mathbf{Y}_t\mid\theta_0)$ implies
$Z_t^T(\theta-\theta_0) \ge\frac{1}2t\mu\llvert  \theta-\theta_0\rrvert  ^2$, which
in turn implies $\llvert  Z_t\rrvert  \ge\frac{1}2t\mu\llvert  \theta-\theta_0\rrvert  $. This is
equivalent to the first statement.

Let us then prove the latter statement. Now $\llvert  Z_t\rrvert  t^{-1/2}=\llvert  A\rrvert  ^t\ge
a$ implies that $\llvert  Z_t^{(k)}\rrvert  \ge t^{1/2}a/\sqrt{n}$ holds for at
least one component $k\in\{1,\ldots,n\}$. But as each $Z_t^{(k)}$ is
a martingale satisfying $Z_0^{(k)}=0$ and
$\llvert  Z_{k+1}^{(k)}-Z_{k}^{(k)}\rrvert  \le M$, Theorem~\ref{thmazuma} yields
\[
\Pr \bigl\{ \bigl\llvert Z_t^{(k)} \bigr\rrvert \ge
t^{1/2}a/\sqrt n \bigr\} \le2\exp \biggl(-\frac{ta^2}{2ntM^2} \biggr)
\]
for all $k\in\{1,\ldots,n\}$. Summing these probabilities over $k$
so as to give an upper bound on the probability that at least one
component is over the limit gives the statement.
\end{pf}
%

\begin{lemma}\label{lemL3}
Suppose that \textup{N1--N3} and \textup{L1} hold. Then, \textup{L3} holds almost surely.
\end{lemma}

\begin{pf}
For any sufficiently small $\varepsilon>0$, N2 implies that the set
$V=B(\theta_0,\varepsilon)$ is a subset of $\mathsfTheta $.
Lemma~\ref{lemoutprob} applied to this set implies that
$\Theta^*_t$ converges fast in probability to $\theta_0$, that is,
the probability $\Pr\{\Theta^*_t\notin B(\theta_0,\varepsilon)\}$
sums to a
finite value over all $t$. This implies that
$\Theta^*_t\asto\theta_0$.
\end{pf}


\begin{theorem}[(Asymptotic normality)]\label{thmasymnorm}
Suppose that \textup{N1--N5} hold and let \textup{L1--L3} hold for some $\mu>0$,
$C>0$, and $U\subset\mathsfTheta $. Then, the following
conditions surely hold when $t$ increases within indices satisfying
$\lambda_t\ge t\mu$:
\begin{enumerate}[3.]
\item The posterior density of the scaled variable $\Phi_t =
t^{1/2}(\Theta- \Theta^*_t)$ satisfies
\[
\int \bigl\llvert p_t(\phi_t\mid\Theta\in U) - N \bigl(
\phi_t; 0, B_t^{-1} \bigr) \bigr\rrvert \,\mathrm{d}\phi
_t \to0,
\]
where $N(\cdots)$ denotes a normal density with given mean and
covariance and $B_t = -t^{-1}\nabla^2_\theta \log
p(\mathbf{Y}_t\mid\theta_0)$.
\item All moments as well as the entropy of $p_t(\phi_t\mid
\Theta\in U)$ are asymptotically equal to those of
$N(\phi_t;0,B_t^{-1})$, that is, the difference converges to zero.
\item
Adjusting\vspace*{1pt} for the $t^{1/2}$ scaling factor, this implies in
particular that 
$t\Cov_t(\Theta\mid U) - B_t^{-1}\to0$ and
$t^{3/2}\mathrm{E}_t(\llvert  \Theta-\mathrm{E}_t(\Theta\mid U)\rrvert  ^3\mid
U)\le
c_n\mu^{-3/2}$ for sufficiently large $t$ for some constant
$c_n$, and so (assuming that $U$ is bounded and convex),
Lemma~\ref{lemiapprox} yields
\[
\sup_{x\in\mathsf{X}} \biggl\llvert t\I_t(
\Theta;Y_x\mid U)-\frac{1}2B_t^{-1}
\odot I_x(\theta_0) \biggr\rrvert \to0.
\]
\end{enumerate}
\end{theorem}

\begin{pf}
The scaled variable $\Phi_t$ takes values in the set $V_t:=
\{ \phi_t\in\mathbb{R}^n\dvt  \Theta^*_t+t^{-1/2}\phi_t\in U \}$. A~Taylor expansion of $\log p(\mathbf{Y}_t\mid\phi_t)$ at $\phi_t =
0$ yields
%
\[
\frac{p_t(\phi_t)}{p_t(\phi_t=0)} = \exp \bigl(\pm\varepsilon (r) \bigr)\exp \biggl( -
\frac{1}2\phi_t^T B_t
\phi_t \pm\frac{1}2\varepsilon(r)\llvert \phi_t
\rrvert ^2 \biggr)
\]
for all $\phi_t$ satisfying $\Theta^*_t+t^{-1/2}\phi_t\in
B(\theta_0,r)$, where
\[
\varepsilon(r) = \sup_{x,y_x,\theta\in B(\theta_0,r)} \max \biggl\{ \biggl\llvert \log
\frac{p(\theta)}{p(\theta')} \biggr\rrvert, \bigl\llvert \nabla^2_\theta
\log p(y_x\mid\theta) - \nabla^2_\theta\log p
\bigl(y_x\mid\theta' \bigr) \bigr\rrvert \biggr\}.
\]
Denoting $r_t = t^{1/4}$, we have $S_t:= B(0,r_t)\subset V_t$ for
sufficiently large $t$ and $\varepsilon_t =
\varepsilon(r_tt^{-1/2}+\llvert  \Theta^*_t-\theta_0\rrvert  )\to0$. It follows
\[
p_t(\phi_t)\propto f_t(\phi_t)
:= \underbrace{\exp \bigl(-\tfrac{1}2\phi_t^TB_t
\phi_t \bigr)}_{=:N_t(\phi_t)}g_t(\phi_t)
\]
for all $\phi_t\in V_t$, where $g_t(\phi) =
\exp(\pm\varepsilon_t\pm\frac{1}2\varepsilon_t^{1/2})\to1$ for
$\phi\in S_t$. As
$N_t(\phi)$ is uniformly bounded and $S_t\to\mathbb{R}^n$, it
follows $[\phi\in V_t]f_t(\phi)-N_t(\phi)\to0$ for all
$\phi\in\mathbb{R}^n$. Furthermore, as $N_t(\phi_t)\le
\exp(-\frac{1}2\mu\llvert  \phi\rrvert  ^2)$ and $g_t(\phi) =
\exp(\pm C\pm\frac{1}4\mu\llvert  \phi\rrvert  ^2)$ for all $\phi\in V_t$, it follows
\[
\int[\phi\in V_t]f_t(\phi)\llvert \phi\rrvert
^k \le \int \exp \biggl(C-\frac{1}4\mu\llvert \phi\rrvert
^2 \biggr)\llvert \phi\rrvert ^k <\infty, \qquad \int
N_t(\phi)\llvert \phi\rrvert ^k <\infty
\]
for all $k\ge0$, and so Lebesgue's dominated convergence theorem
implies that
\[
\int \bigl\llvert [\phi\in V_t]f_t(\phi)u(\phi) -
N_t( \phi)u(\phi) \bigr\rrvert \,\mathrm{d}\phi\to0
\]
for any function $\llvert  u(\phi)\rrvert  \le\llvert  \phi\rrvert  ^k$. This implies that all
moments of $[\phi\in V_t]f_t(\phi)$ are asymptotically equal to
those of $N_t(\phi)$. As the eigenvalues of $B_t$ are between $\mu$
and $M$, the normalization constant $Z:=\int N_t(\phi)\,\mathrm{d}\phi$ is
within the constant range $[(2\uppi/M)^{n/2},(2\uppi/\mu)^{n/2}]$, and
it follows that the moments of the normalized densities
$p_t(\phi_t)$ and $N(\phi_t;0,B_t^{-1})$ are also asymptotically
equal. Similarly, as $f_t(\phi)\log f_t(\phi) - N_t(\phi)\log
N_t(\phi)\to0$, where the log-factors can be bounded by
polynomials of $\llvert  \phi\rrvert  $, it follows that the entropies of
$p_t(\phi_t)$ and $N(\phi_t;0,B_t^{-1})$ are asymptotically equal.
(Note that the entropy of a density $p(x)=f(x)/Z$ can be calculated
as $-(\int f\log f)/Z + \log(Z)$.)
\end{pf}


\begin{lemma}\label{lemasymhess}
Suppose that \textup{N1} and \textup{N3} hold. Then, conditioned on $\theta_0$ as the
true parameter value, $\mathrm{E}(-\nabla^2_\theta \log p(Y_x\mid
\theta_0))
= I_x(\theta_0)$ for all $x\in\mathsf{X}$, and
\[
B_t - \frac{\sum_{k=1}^t I_{X_t}(\theta_0)}{t}\astoi0,
\]
where $B_t = -t^{-1}\nabla^2_\theta \log
p(\mathbf{Y}_t\mid\theta_0)$.
\end{lemma}

\begin{pf}
\begin{eqnarray*}
&& \mathrm{E}\bigl(-\nabla^2_\theta \log p(Y_x\mid\theta_0)\mid\Theta=
\theta_0\bigr)
\\
&&\quad = \int p(y_x\mid\theta_0) \biggl\{
\biggl[\frac{\nabla_\theta p(y_x\mid\theta_0)}{p(y_x\mid\theta
_0)} \biggr] \biggl[\frac{\nabla_\theta p(y_x\mid\theta_0)}{p(y_x\mid\theta
_0)} \biggr]^T
- \frac{\nabla^2_\theta p(y_x\mid\theta_0)}{p(y_x\mid\theta_0)} \biggr\}\,\mathrm{d}y_x
\\
&&\quad = I_x(\theta_0) - \int\nabla^2_\theta
p(y_x\mid\theta_0)\,\mathrm{d}y_x
\\
&&\quad = I_x(\theta_0) - \nabla_\theta\int
\nabla_\theta p(y_x\mid \theta_0)\,\mathrm{d}y_x
\\
&&\quad = I_x(\theta_0) - \nabla^2_\theta
\int p(y_x\mid\theta_0)\,\mathrm{d}y_x =
I_x(\theta_0),
\end{eqnarray*}
where the interchange of the order of integration and
differentiation is justified by Lebesgue's dominated convergence
theorem for the $\mathrm{d}y_x$-integrable dominating functions $f_x(y_x)$
and $g_x(y_x)$ given by
%
%
\begin{eqnarray*}
\bigl\llvert \nabla^2_\theta p(y_x\mid\theta)
\bigr\rrvert &=& p(y_x\mid\theta) \bigl\llvert \nabla_\theta
\log p(y_x\mid\theta) \nabla_\theta \log p(y_x
\mid\theta)^T + \nabla^2_\theta \log
p(y_x\mid\theta) \bigr\rrvert
\\
&\le& p(y_x\mid\theta_0)\exp \bigl(M\llvert \theta-
\theta_0\rrvert \bigr)\cdot \bigl(M^2+M \bigr)
\\
&\le& p(y_x\mid\theta_0)\exp(MD)\cdot
\bigl(M^2+M \bigr) =: f_x(y_x)
\end{eqnarray*}
and
%
%
\begin{eqnarray*}
\bigl\llvert \nabla_\theta p(y_x\mid\theta) \bigr\rrvert
&=& p(y_x\mid\theta) \bigl\llvert \nabla_\theta \log
p(y_x\mid\theta) \bigr\rrvert
\\
&\le& p(y_x\mid\theta_0)\exp(MD)\cdot M =:
g_x(y_x).
\end{eqnarray*}
Thus, denoting $Z_k = -\nabla^2_\theta \log p(Y_{x_k}\mid\theta_0)
- I_{X_k}(\theta_0)$, given $\Theta=\theta_0$, the sequence
$Z_1+\cdots+Z_k$ of partial sums is a martingale and satisfies
$\mathrm{E}(\llvert  Z_k\rrvert  ^2)\le(M+M)^2 <\infty$ for all $k$, and so
Theorem~\ref{thmmartingaleSLLN} implies that
$(Z_1+\cdots+Z_t)/t\asto0$, which is the statement.
\end{pf}

%
\begin{corollary}\label{corvarlower}
Suppose that \textup{N1--N5} hold. Then, for all $\mu>0$, almost surely
$t\Sigma_t>(B_t+\mu I)^{-1}$ (meaning that the difference is
positive definite) for all sufficiently large $t$, where $B_t:=
-t^{-1}\nabla^2_\theta \log p(\mathbf{Y}_t\mid\theta_0)$. In
particular, $\tr(t\Sigma_t) \ge(2\mu)^{-1}$ and $\det(t\Sigma
_t)\ge
(2\mu)^{-1}(2M)^{-(n-1)}$ for all sufficiently large $t$ satisfying
$\min\lambda_{B_t}\le\mu\le M$, where $\min\lambda_{B_t}$ denotes
the smallest eigenvalue of $B_t$.
\end{corollary}
\begin{pf}
Let $\mu>0$ be arbitrary and define an augmented observation model
$Y'_x:= (Y_x,Z)$, where $Z\sim N(\Theta,\mu^{-1}I)$ is independent
(given $\theta$) from $Y_x$. Let $U$ be a neighborhood of
$\theta_0$ satisfying L1 and L2 as well as L3 almost surely. If we
choose the auxiliary component $z_t$ so as to obtain
$t^{-1}\sum_{k=1}^t z_k = \mathrm{E}(\Theta\mid\mathbf{y}_t)$ for
each $t$,
then L3 remains satisfied given the augmented data and we also
obtain $\Sigma_t>\Sigma'_t$, because the augmented data will
strictly decrease the square error from the original mean, and
moving to the new mean can only further reduce this error. The
normalized Hessian at $\theta_0$ for the augmented data is $B'_t =
B_t + \mu I$, and so, due to Lemma~\ref{lemasymhess},
$\min\lambda_{B'_t}\ge\mu/2$ for all sufficiently large $t$
(although we have fiddled with the $z_k$ values,
Lemma~\ref{lemasymhess} still applies as it does not depend on
these values). Thus, Theorem~\ref{thmasymnorm}(3) implies that
$t\Cov(\Theta\mid\mathbf{y}'_t, U) - (B'_t)^{-1}\to0$ (note that
Theorem~\ref{thmasymnorm} is a sure result and hence applies even
with our fiddled $z_k$ values). Since $\Pr_t\{\Theta\in U^c\}$
decays exponentially in the augmented model, it follows that also
$t\Sigma_t'-(B'_t)^{-1}\to0$. As the eigenvalues of $B'_t$ are
within the range $[\mu/2,M+\mu/2]$, the matrix inverse behaves
nicely and we obtain $(t\Sigma_t')^{-1} - B_t'\to0$, which implies
$(t\Sigma_t')^{-1} - B_t' < \varepsilon I$ for all sufficiently large $t$
for any $\varepsilon>0$. It follows $t\Sigma_t>t\Sigma_t'>
(B_t+(\mu+\varepsilon)I)^{-1}$ for all sufficiently large $t$.
\end{pf}


\section{Asymptotic optimality}\label{secoptimality}

In this section, we assume that:
\begin{enumerate}[O4.]
\item[O1.] \textup{C1--C4} hold globally.
\item[O2.] Some neighborhood $U_0$ of $\theta_0\in\mathsfTheta $ is
homeomorphic to a subset of $\mathbb{R}^n$ that satisfies \mbox{\textup{N1--N5}}.
\item[O3.] There exists placements $x_1,\ldots,x_m\in\mathsf{X}$ and
nonnegative weights \mbox{$\alpha_1+\cdots+\alpha_m=1$} such that
$\sum_{j=1}^m \alpha_j I_{x_j}(\theta_0)$ is positive definite.
\item[O4.] The placements $X_t$ satisfy
\[
R_t:= \frac{
\I_t(\Theta;Y_{X_{t+1}})
}{
\sup_{x\in\mathsf{X}}\I_t(\Theta;Y_x)
}\leadsto1.
\]
(See Section~\ref{secasymprop} below for the definition of
``$\leadsto$''.)
\end{enumerate}

First, let us say a few words about the main difficulty related to the
adaptivity of the placements, namely the complications caused by any
secondary modes in the posterior distribution. This issue is
discussed by Paninski \cite{paninski2005} in the context of consistency, but
it seems that even after consistency has been established, the issue
cannot be ignored.

The information maximization strategy decreases the relative weights
of any secondary modes only at a rate approximately proportional to
$1/t$ \cite{paninski2005}. Therefore, any secondary mode may have a
contribution proportional to $1/t$ to all moments of the posterior
distribution. This means that only the first order moments of the
approximating normal distribution remain asymptotically accurate, even
though its total variation distance from the posterior does tend to
zero. In particular, the inverse Hessian of the likelihood generally
\emph{does not} give an asymptotically accurate approximation of the
global posterior covariance. (In fact, the global posterior
covariance may be undefined as $\mathsfTheta $ need not have a
global Euclidean structure.)

For this reason, the asymptotic approximation to the expected
information gain $\I_t(\Theta;Y_x\mid U)$ given by
Theorem~\ref{thmasymnorm}(3) only applies within a sufficiently small
neighborhood $U$ of the true parameter value, where the posterior can
be shown to be asymptotically unimodal. Nonetheless, even though the
local and global moments are not in good agreement asymptotically, it
turns out that $\I_t(\Theta;Y_{X_{t+1}}\mid U)$ is in fact in good
agreement with $\I_t(\Theta;Y_{X_{t+1}})$ on ``most trials''. Indeed,
as the relative weights of any secondary modes typically decay at an
exponential rate with the number of trials whose placements can
distinguish between them, it follows that the placements of only a
decreasing fraction of trials can be significantly affected by the
secondary modes.

To formalize this intuition, we will first develop a theory for
measuring asymptotic proportions.

\subsection{Asymptotic proportions}\label{secasymprop}

%
\begin{definition}
To measure subsets $K\subset\mathbb{N}$, we use the \emph{proportion
measures}
\[
\rho(K)=\lim_{n\to\infty}\rho_{1,n}(K),\qquad
\rho_{a,b}(K)=\frac{\llvert  K\cap[a,b[ \rrvert  }{b-a},
\]
where $\llvert  \cdot\rrvert  $ indicates the cardinality of a set. (Note that
although $\rho_{a,b}$ is a measure in the measure-theoretic sense
for any $a,b\in\mathbb{N}$, the limit $\rho$ is only a
\emph{finitely additive} measure.) When we say ``for almost every
$n\in\mathbb{N}$'', we mean that the set where the statement does
not hold is a null set w.r.t. $\rho$. We use the notation
$x_k\leadsto x$ to mean that there exists a subset
$K\subset\mathbb{N}$ with $\rho(K)=1$ such that $[k\in K](x_k-x)\to
0$. We also define
%
%
\begin{eqnarray*}
\limsup_{k\leadsto\infty}x_k&:=& \inf\{ x\in
\mathbb{R}\dvt x_k\le x\mbox{ for a.e. }k\in\mathbb{N} \},
\\
\liminf_{k\leadsto\infty}x_k&:=& \sup\{ x\in
\mathbb{R}\dvt x_k\ge x\mbox{ for a.e. }k\in\mathbb{N} \},
\end{eqnarray*}
and when both equal $x$, we write $\lim_{k\leadsto\infty}x_k = x$.
\end{definition}

%
\begin{lemma}\label{lempropconj}
Suppose that for all $j\in\mathbb{N}$, the proposition $P_k^j$ holds
for a.e. $k\in\mathbb{N}$. Then there exists an increasing
sequence $j(k)\to\infty$ such that $P_k^1\wedge\cdots\wedge
P_k^{j(k)}$ holds for a.e. $k\in\mathbb{N}$.
\end{lemma}

\begin{pf}
For all $j\in\mathbb{N}$, $Q_k^j:= P_k^1\wedge\cdots\wedge P_k^j$
holds for a.e. $k\in\mathbb{N}$. Thus, for all $j\in\mathbb{N}$,
\[
f_j(k):= \inf_{k'\ge k}\frac{\sum_{i=1}^{k'}Q_i^j}{k'}
\]
is increasing in $k$ and tends to one as $k\to\infty$. Choosing
\[
j(k) = \max \bigl\{ j'\in\mathbb{N}\dvt  f_{j'}(k)
\ge1-1/j' \bigr\}
\]
yields the statement.
\end{pf}
%

%
\begin{lemma}\label{lemproplemma}
If $x_k$ is a bounded sequence, then the following are equivalent:
\begin{enumerate}[4.]
\item$x_k\leadsto x$,
\item$\llvert  x_k-x\rrvert  <\varepsilon$ for a.e. $k\in\mathbb{N}$ for all
$\varepsilon>0$,
\item$\lim_{k\leadsto\infty}x_k = x$,
\item$ \frac{1}t\sum_{k=1}^t\llvert  x_k-x\rrvert  \to0$.
\end{enumerate}
If $x_k$ is not bounded, then 1--3 are equivalent and implied by 4.
\end{lemma}
\begin{pf}
All implications are fairly obvious. As an example, ``2
$\Rightarrow$ 1'' follows from Lemma~\ref{lempropconj} applied to
$P^j_k = [ \llvert  x_k-x\rrvert  <1/j ]$.
\end{pf}

%
\begin{lemma}\label{lemsubharmonic}
Let $x_k$ be a nonnegative sequence. If $\sum_{k=1}^\infty x_k <
\infty$, then for any $\varepsilon>0$, the inequality $x_k <
\varepsilon/k$ holds
true for almost every $k\in\mathbb{N}$ (which implies $k\cdot
x_k\leadsto0$).
\end{lemma}
\begin{pf}
Assume the contrary: for some $\varepsilon>0$ there exists a set
$K\subset\mathbb{N}$ such that $x_k\ge\varepsilon/k$ for all $k\in
K$ and
for some $c>0$, $\rho_{1,k}(K)> c$ for arbitrarily large $k$. As
$\rho_{1,n+1}(K) - \rho_{k,n+k}(K)\le2k/n \to0$ as $n\to\infty$
for all $k$, we can recursively find an increasing sequence of
indices $k_1=1$, $k_{i+1}\ge2k_i$, such that
$\rho_{k_i,k_{i+1}}(K)\ge c$ for all $i$. This yields
%
%
\begin{eqnarray*}
\sum_{k=1}^\infty x_k \ge\sum
_{i=1}^\infty c(k_{i+1}-k_i)
\frac \varepsilon{k_i} \ge\sum_{i=1}^\infty
c(2k_i-k_i)\frac\varepsilon{k_i} = \infty,
\end{eqnarray*}
which contradicts the assumption.
\end{pf}

%
\begin{lemma}
Suppose that a sequence of random variables $X_k\dvtx \Omega\to[-M,M]$
satisfies $X_k\leadsto X$ almost surely. Then, $\mathrm
{E}(\llvert  X_k-X\rrvert  )\leadsto0$.
\end{lemma}
\begin{pf}
By Lemma~\ref{lemproplemma}(4) and the dominated convergence
theorem,
\[
\frac{1}t\sum_{k=1}^t \mathrm{E}
\bigl(\llvert X_k-X\rrvert \bigr) = \mathrm{E} \Biggl(
\frac{1}t\sum_{k=1}^t\llvert
X_k-X\rrvert \Biggr) \to \mathrm{E} \Biggl(\lim_{t\to\infty}
\frac{1}t\sum_{k=1}^t\llvert
X_k-X\rrvert \Biggr)=0.
\]
\upqed
\end{pf}

%
\begin{corollary}
Suppose that the event $A_k$ happens for a.e. $k\in\mathbb{N}$ a.s.
Then, $\Pr\{A_k\}\leadsto1$.
\end{corollary}

%
\begin{definition}
We use the notation $X_k\prhoto X$ to mean that there exists a
subset $K\subset\mathbb{N}$ with $\rho(K)=1$ such that $[k\in
K](X_k-X)\pto0$.
\end{definition}

%
\begin{lemma}\label{lemprhotoequiv}
$X_k\prhoto X$ if and only if
$\Pr\{\llvert  X_k-X\rrvert  \ge\varepsilon\}\leadsto0$ for all $\varepsilon>0$.
\end{lemma}
\begin{pf}
The ``only if'' direction is obvious. We will prove the ``if''
direction.

By definition, we have $\Pr\{\llvert  X_k-X\rrvert  \ge1/j\}\le1/j$ for a.e. $k\in\mathbb{N}$ for all $j\in\mathbb{N}$. Lemma~\ref{lempropconj}
then implies that there exists an increasing sequence
$j(k)\to\infty$ such that
\[
\Pr \bigl\{\llvert X_k-X\rrvert \ge1/j(k) \bigr\}\le1/j(k)\to0
\]
for a.e. $k\in\mathbb{N}$.
\end{pf}

%
\begin{lemma}\label{lemasprhoto}
Suppose that a sequence of random variables $X_k$ satisfies
$X_k\leadsto X$ almost surely. Then, $X_k\prhoto X$.
\end{lemma}
\begin{pf}
Let $\varepsilon>0$ be arbitrary. Denoting
\[
Y_t = \frac{1}t\sum_{k=1}^t
\bigl[ \llvert X_k - X\rrvert \ge\varepsilon \bigr],
\]
$X_k\leadsto X$ implies that $Y_t\to0$. As $Y_t$ is bounded, the
dominated convergence theorem implies
\[
0 = \mathrm{E} \Bigl(\lim_{t\to\infty} Y_t \Bigr) = \lim
_{t\to
\infty} \mathrm{E}(Y_t) = \lim_{t\to\infty}
\frac{1}t\sum_{k=1}^t\Pr \bigl\{
\llvert X_k-X\rrvert \ge\varepsilon \bigr\}
\]
and so Lemma~\ref{lemproplemma}(4) yields
$\Pr\{\llvert  X_k-X\rrvert  \ge\varepsilon\}\leadsto0$. Now Lemma~\ref{lemprhotoequiv}
implies the statement.
\end{pf}

\subsection{Asymptotic D-optimality}\label{secd-optimality}

In this section, we show that the greedy information maximization
strategy satisfies asymptotically a condition known as
D-optimality. This condition is defined as maximality of the
determinant of the Fisher information matrix of the experiment at the
true parameter value $\theta_0$. The D-optimality criterion is
special among all functionals of the information matrix (such as the
trace, minimum eigenvalue, etc.) in that it is insensitive to linear
or affine transformations of the parameter space $\mathsfTheta $.
Furthermore, in the asymptotically normal models that we are
interested in, it yields a (local) approximation of the posterior
entropy, which is the utility function commonly used in adaptive
estimation settings. We will make use of this fact in the next section
to derive an asymptotic expression of the posterior entropy.

%
\begin{lemma}\label{leminfolower}
For almost any $\theta_0\in\mathsfTheta $ satisfying \textup{O1--O3}, there
exists a constant $c$ such that for all $\mu>0$, given $\theta_0$ as
the true parameter value, almost surely $\I_t(\Theta;Y_{X_{t+1}})\ge
c(t\mu)^{-1}$ for all sufficiently large $t$ satisfying
$\lambda_t\le t\mu$, where $\lambda_t$ denotes the smallest
eigenvalue of $-\nabla^2_\theta \log p(\mathbf{Y}_t\mid\theta_0)$.
\end{lemma}

\begin{pf}
Denoting\vspace*{1pt} $I:= \sum_{j=1}^m\alpha_j I_j$, where $\alpha_j$ and $I_j:= I_{x_j}(\theta_0)$ are given by O3, the smallest eigenvalue
$\min\lambda_I$ is positive.

Suppose that $U_0$ has diameter $D$ and let $C_{M,D}$ be the constant
of Lemma~\ref{lemiapprox} applied to $U_0$ as the parameter space.
The same constant also applies to any subset $U =
B(\theta_0,\delta/2)\subset U_0$ with diameter $\delta\le D$ and as
the posteriors are strongly consistent in $U$, too,
Lemma~\ref{lemscimpl} implies that $\mathrm{E}_t(\Theta\mid
U)\asto\theta_0$. Thus, N3 and N4 imply that $\llvert  I_x(\mathrm
{E}_t(\Theta\mid
U))-I_x(\theta_0)\rrvert  <\delta$ for all\vadjust{\goodbreak} $x$ for all sufficiently large
$t$. We obtain
%
%
\begin{eqnarray*}
\I_t(Y_x;\Theta\mid U) &\ge&\frac{1}2
\Cov_t(\Theta\mid U)\odot I_x \bigl(\mathrm{E}_t(
\Theta\mid U) \bigr) - C_{M,D}\mathrm{E}_t \bigl( \bigl
\llvert \Theta - \mathrm{E}_t(\Theta\mid U) \bigr\rrvert
^3 \mid U \bigr)
\\
&\ge& \frac{1}2 \Cov_t(\Theta\mid U)\odot I_x
\bigl(\mathrm{E}_t(\Theta\mid U) \bigr) - C_{M,D}
\mathrm{E}_t \bigl(\delta \bigl\llvert \Theta-\mathrm{E}_t(
\Theta\mid U) \bigr\rrvert ^2\mid U \bigr)
\\
&=& \frac{1}2 \tr \bigl(\Cov_t(\Theta\mid U)I_x
\bigl(\mathrm{E}_t(\Theta\mid U) \bigr) \bigr) - C_{M,D}
\delta\tr \bigl(\Cov_t(\Theta\mid U) \bigr)
\\
&\ge& \frac{1}2 \tr \bigl(\Cov_t(\Theta\mid
U)I_x(\theta_0) \bigr) - \biggl(C_{M,D}+
\frac{1}2 \biggr)\delta\tr \bigl(\Cov_t(\Theta\mid U) \bigr)
\\
&\ge& \frac{1}2 \max_{j=1,\ldots,m}\tr \bigl(
\Cov_t(\Theta\mid U)I_j \bigr) - \biggl(C_{M,D}+
\frac{1}2 \biggr)\delta\tr \bigl(\Cov_t(\Theta\mid U) \bigr)
\\
&\ge& \frac{1}2\tr \bigl(\Cov_t(\Theta\mid U)I \bigr) -
\biggl(C_{M,D}+\frac{1}2 \biggr)\delta\tr \bigl(
\Cov_t(\Theta\mid U) \bigr)
\\
&\ge& \frac{1}2 \tr \bigl(\Cov_t(\Theta\mid U) \bigr) \min
\lambda_I - \biggl(C_{M,D}+\frac{1}2 \biggr)\delta
\tr \bigl(\Cov_t(\Theta\mid U) \bigr)
\\
&=& \biggl(\frac{\min\lambda_I}{2}- \biggl(C_{M,D}+\frac{1}2
\biggr) \delta \biggr)\tr \bigl(\Cov_t(\Theta\mid U) \bigr) =: c \tr
\bigl( \Cov_t(\Theta\mid U) \bigr)
\end{eqnarray*}
for some $x\in\mathsf{X}$ (fourth inequality) for all sufficiently
large $t$ (third inequality), where we have used the fact that
$\tr(A)\min\lambda_B\le\tr(AB)\le\tr(A)\max\lambda_B$ (sixth and
third inequalities). Let us then choose $\delta<
\min\lambda_I/(2C_{M,D}+1)$ so that $c$ as defined above is positive.
Now, the inequality $\I_t(\Theta;Y_x)\ge p_t(U)\I_t(\Theta;Y_x\mid
U)$, which follows from the chain rule of mutual information (cf. the proof of the next lemma), and \textup{C4}${}+{}$Corollary~\ref{corvarlower}
imply
%
%
\begin{eqnarray*}
\I_t(\Theta;Y_{t+1}) &\ge& \gamma\sup_{x\in\mathsf{X}}
\I_t(\Theta;Y_x) \ge \gamma\sup_{x\in\mathsf{X}}p_t(U)
\I_t(\Theta;Y_x\mid U)
\\
&\ge&\gamma p_t(U)c\tr \bigl(\Cov_t(\Theta\mid U) \bigr)
\ge \gamma p_t(U)c(2t\mu)^{-1}.
\end{eqnarray*}
As Lemma~\ref{lemscimpl} yields $p_t(U)\asto1$, the statement
follows.
\end{pf}

%
\begin{lemma}\label{leminfolocal}
For almost any $\theta_0\in\mathsfTheta $ satisfying \textup{O1--O3}, there
exists a neighborhood $U\subset U_0$ of~$\theta_0$ such that
conditioned on $\theta_0$ as the true parameter value, almost
surely,
\[
Q_t:= \frac{\I_t(\Theta;Y_{X_{t+1}}\mid U)}{\I_t(\Theta;Y_{X_{t+1}})}\leadsto1.
\]
\end{lemma}

\begin{pf}
By Lemmas~\ref{lemfinitesum},~\ref{leminfoconv} and
\ref{lemsubharmonic}, almost surely, the convergences
%
%
\begin{eqnarray*}
\I_t(\Theta;Y_{X_{t+1}}\mid U) &\to&0,
\\
t \I_t \bigl([\Theta\in U]; Y_{X_{t+1}} \bigr) &\leadsto& 0
\end{eqnarray*}
hold for all neighborhoods $U$ in a countable basis of the compact
metrizable space $\mathsfTheta $. It follows that the same is
true conditioned on almost any $\theta_0\in\mathsfTheta $ as the
true parameter value. Thus, given almost any
$\theta_0\in\mathsfTheta $, we can pick a neighborhood $U\subset
U_0$ of $\theta_0$ from the countable basis such that the above
convergences almost surely hold.

Lemma~\ref{leminfolower} (applied to $\mu= M$) almost surely
yields
\[
\I_t(\Theta;Y_{t+1})\ge c(M t)^{-1}=:c_1t^{-1}
\]
for all sufficiently large $t$, where we denote $Y_{t+1} = Y_{X_{t+1}}$.
Condition~C4${}+{}$Lemma~\ref{lemscimpl} yields
\[
\I_t(\Theta;Y_{t+1})\ge\gamma\sup_{x\in\mathsf{X}}
\I_t(\Theta;Y_x) \ge\gamma c p_t
\bigl(U^c \bigr) =: c_2 p_t
\bigl(U^c \bigr)
\]
for all sufficiently large $t$, and the chain rule of mutual
information yields
\[
\I_t(\Theta;Y_{t+1}) = \I_t \bigl([\Theta\in
U]; Y_{t+1} \bigr) + p_t(U)\I_t(
\Theta;Y_{t+1} \mid U) + p_t \bigl(U^c \bigr)
\I_t \bigl( \Theta;Y_{t+1}\mid U^c \bigr).
\]
Thus, almost surely,
\[
\frac{\I_t(\Theta;Y_{t+1}\mid U)}{\I_t(\Theta;Y_{t+1})} = \underbrace{\frac{1}{p_t(U)}}_{\to1} \biggl[1 -
\frac{
\overbrace{p_t(U^c)}^{\le\I_t(\Theta;Y_{t+1})/c_2}
\overbrace{\I_t(\Theta;Y_{t+1}\mid U^c)}^{\to0} +
\overbrace{t\I_t([\Theta\in U]; Y_{t+1})}^{\leadsto0}t^{-1}
}{
\underbrace{\I_t(\Theta;Y_{t+1})}_{
\ge c_1t^{-1}
}
} \biggr]\leadsto1.
\]\upqed
\end{pf}

%
\begin{corollary}\label{corinfohess}
Conditioned on almost any $\theta_0$ satisfying \textup{O1--O4}, the sequence
\[
D_t:= \sup_{x\in\mathsf{X}}B_t^{-1}
\odot I_x(\theta_0) - B_t^{-1}
\odot \II_{X_{t+1}}(\theta_0)
\]
satisfies $[\min\lambda_{B_t}\ge\mu]D_t\leadsto0$ a.s. for any given
$\mu>0$, where $\min\lambda_{B_t}$ denotes the smallest eigenvalue
of $B_t:= -t^{-1}\nabla^2_\theta \log
p(\mathbf{Y}_t\mid\theta_0)$.
\end{corollary}
\begin{pf}
Let us first shrink the neighborhood $U_0$ of $\theta_0$ as
necessary to make its diameter smaller than the constant
$\delta_{\mu,C}$ given by Lemma~\ref{lemL12}. Then, let $U\subset
U_0$ be the neighborhood of $\theta_0$ given by
Lemma~\ref{leminfolocal}. By Theorem~\ref{thmasymnorm}(3), there
now exist random sequences $E_t\to0$ and $E_t'\to0$ such that
conditioned on $\theta_0$ as the true value,
%
%
\begin{eqnarray*}
\frac{1}2 \sup_{x\in\mathsf{X}} B_t^{-1}
\odot I_x(\theta_0) &=& \sup_{x\in\mathsf{X}} t
\I_t(\Theta;Y_x\mid U) + E_t,
\\
\frac{1}2 B_t^{-1}\odot \II_{X_{t+1}}(
\theta_0) &=& t\I_t(\Theta;Y_{X_{t+1}}\mid U) +
E_t'
\end{eqnarray*}
whenever $\min\lambda_{B_t}\ge\mu$. For these $t$, it follows
%
%
\begin{eqnarray*}
\frac{1}2D_t &=& \biggl(\frac{1}2 \sup
_{x\in\mathsf{X}} \underbrace{B_t^{-1}\odot
I_x(\theta_0)}_{
= \tr(B_t^{-1}I_x(\theta_0))\le n\mu^{-1}M
} - E_t
\biggr) \biggl(1 - \frac{
\I_t(\Theta;Y_{X_{t+1}}\mid U)
}{
\sup_{x\in\mathsf{X}}\I_t(\Theta;Y_x\mid U)
} \biggr) + E_t-E_t',
\end{eqnarray*}
where Lemma~\ref{leminfolocal} and the inequality
$\I_t(\Theta;Y_x)\ge p_t(U)\I_t(\Theta;Y_x\mid U)$ yield
%
%
\begin{eqnarray*}
\frac{
\I_t(\Theta;Y_{X_{t+1}}\mid U)
}{
\sup_{x\in\mathsf{X}}\I_t(\Theta;Y_x\mid U)
} \ge p_t(U)\frac{
\I_t(\Theta;Y_{X_{t+1}}\mid U)
}{
\sup_{x\in\mathsf{X}}\I_t(\Theta;Y_x)
} = p_t(U)
Q_t R_t \leadsto1,
\end{eqnarray*}
and so $[\min\lambda_{B_t}\ge\mu]D_t\leadsto0$.
\end{pf}

%
\begin{lemma}\label{lemminlambda}
Conditioned on almost any $\theta_0$ satisfying \textup{O1--O3}, there exists
$\mu$ such that $\min\lambda_{B_t}\ge\mu$ for infinitely many
$t\in\mathbb{N}$, where $\min\lambda_{B_t}$ denotes the smallest
eigenvalue of $B_t =
-t^{-1}\nabla^2_\theta p(\mathbf{Y}_t\mid\theta_0)$.
\end{lemma}
\begin{pf}
Let $\mu>0$ be arbitrary. Lemma~\ref{leminfolower} almost surely
yields $\I_{t-1}(\Theta;Y_{X_t})\ge c(t\mu)^{-1}$ for all sufficiently
large $t$ satisfying $\min\lambda_{B_t}<\mu$ and
Lemma~\ref{leminfolocal} implies that $\I_{t-1}(\Theta;Y_{X_t}\mid
U_0)\ge c(t\mu)^{-1}$ for a.e. $t$ satisfying
$\min\lambda_{B_t}\le\mu$. Let then $K_\mu:=
\{ t\in\mathbb{N}\dvt \min\lambda_{B_t}\ge\mu \}$ and suppose that
$\rho(K_\mu)=0$. Then, $\rho_j:= \rho_{2^j,2^{j+1}}(K_\mu)\to0$,
and then
exists $j_0$ such that $\rho_j\le1/2$ for all $j\ge j_0$. It
follows
\[
\sum_{t=1}^{2^{j_1}-1} \I_{t-1}(
\Theta;Y_{X_t}\mid U_0) \ge \frac{c}\mu\sum
_{t=1}^{2^{j_1}-1} [t\notin K_\mu]
\frac{1}t \ge \frac{c}\mu\sum_{j=j_0}^{j_1-1}
\sum_{t=2^j(1+\rho
_j)}^{2^{j+1}-1}\frac{1}t \ge
\frac{c}\mu(j_1-j_0)\log\frac{2}{3/2},
\]
and so
\[
\sum_{k=1}^t \I_{k-1}(
\Theta;Y_{X_k}\mid U_0) \ge \biggl(\frac{c}\mu\log
\frac{4}3 \biggr)\log_2(t-1) - c_{c,\mu}
\]
for all $t = 2^j$, $j\ge j_0$. Since $\mu$ was arbitrary, this
implies that the sum grows asymptotically superlogarithmically if
$\rho(K_\mu)=0$ holds for all $\mu>0$. If this event has positive
probability among all $\theta_0\in U_0$, then also
\[
\I(\Theta;\mathbf{Y}_t\mid U_0) = \mathrm{E} \Biggl(
\sum_{k=1}^t \I_{k-1}(
\Theta;Y_{X_k}\mid U_0) \Bigm| U_0 \Biggr)
\]
grows superlogarithmically, contradicting Lemma~\ref{leminfoupper}.
Thus, for almost all $\theta_0\in U_0$ satisfying \textup{O1--O3}, either
$K_\mu$ is not $\rho$-measurable or $\rho(K_\mu)>0$. In either case
$K_\mu$ is infinite.
\end{pf}


\begin{theorem}[(Asymptotic D-optimality, part 1)]\label{thmdopt}
Conditioned on almost any $\theta_0\in\mathsfTheta $ satisfying
\textup{O1--O4}, almost surely,
%
%
\begin{eqnarray*}
B_t:= -t^{-1}\nabla^2_\theta \log p(
\mathbf{Y}_t\mid\theta_0) \to B^*:=
\argmax_{B\in\mathcal{I}} \det(B),
\end{eqnarray*}
where $\mathcal{I}$ is the convex hull of the closure of
$\{I_x(\theta_0)\}_{x\in\mathsf{X}}$. The maximizer $B^*$ is
unique, because the determinant is log-concave on the compact
convex set $\mathcal{I}$. This result is optimal in the sense that
for any strategy of choosing the placements $X_t$ (instead of \textup{O4} and
\textup{C4}), almost surely $\limsup_{t\to\infty}\det(B_t)\le\det(B^*)$.
\end{theorem}

\begin{pf}
The objective function is
\[
f(B) = \cases{ \displaystyle \log\det(B), &\quad $\min\lambda_B >
0$,
\cr
-\infty,&\quad otherwise,}
\]
where $\lambda_B$ denotes the set of eigenvalues of $B$.
Lemma~\ref{lemasymhess} implies that $B_t$ is asymptotically a
convex combination of matrices in the closure of
$\{I_x(\theta_0)\}_{x\in\mathsf{X}}$ and so
$\limsup_{t\to\infty}f(B_t)\le f(B^*)$. Let us then show that this
upper bound is tight.

First, we choose some representation $B^* = \sum_{k=1}^m\alpha_k
I_k$ of the optimum point, where $I_k$ are matrices in the closure
of $\{I_x(\theta_0)\}_{x\in\mathsf{X}}$ and
$\sum_{k=1}^m\alpha_k=1$.

For any symmetric real matrix $B_t$, we have (with slight abuse of
notation)
%
%
\begin{eqnarray*}
\nabla f(B_t) &=& B_t^{-1},
\\
\nabla^2 f(B_t) &=& - \bigl[ \bigl(B_t^{-1}
\bigr)_i \bigl(B_t^{-1} \bigr)_j^T
\bigr]_{i,j}^n,
\\
{ \bigl[\nabla^2 f(B_t) \bigr]} B &=& - \bigl[
\bigl(B_t^{-1} \bigr)_i \bigl(B_t^{-1}
\bigr)_j^T\odot B \bigr]_{i,j}^n =
-B_t^{-1}BB_t^{-1},
\\
B \odot \bigl[\nabla^2 f(B_t) \bigr] B &=& -\tr
\bigl(B_t^{-1}BB_t^{-1}B \bigr),
\end{eqnarray*}
and Taylor's theorem yields
%
%
\begin{eqnarray*}
f(B_{t+1}) &=& f(B_t) + B_t^{-1}
\odot(B_{t+1}-B_t) - \tfrac{1}2\tr
\bigl(B_t^{-1}B'B_t^{-1}B'
\bigr),
\end{eqnarray*}
where $B'$ is between $0$ and $B_{t+1}-B_t$. Denoting $B:=
-\nabla^2_\theta \log(p(Y_{X_{t+1}}\mid\theta_0))$, we obtain
%
%
\begin{eqnarray*}
f(B_{t+1}) - f(B_t) &=& f \biggl(\frac{tB_t+B}{t+1} \biggr)
- f(B_t)
\\
&=& B_t^{-1}\odot\frac{B-B_t}{t+1} -\frac{1}{2}
\underbrace{\tr \bigl(B_t^{-1}B'B_t^{-1}B'
\bigr)}_{\llvert  \cdot\rrvert  \le n4M^2\mu^{-2}(t+1)^{-2}}
\\
&\ge&\frac{1}{t+1} \biggl(B_t^{-1}\odot B-n -
\frac{2nM^2\mu
^{-2}}{t+1} \biggr),
\end{eqnarray*}
for all indices $t$ satisfying $\min\lambda_{B_t}\ge\mu$ for any
$\mu>0$. Denoting by $\lambda_i$ the eigenvalues of $B_t^{-1}B^*$,
Corollary~\ref{corinfohess} now implies that
%
%
\begin{eqnarray*}
B_t^{-1}\odot \II_{X_{t+1}}(\theta_0) +
D_t &=& \sup_{x\in\mathsf{X}}B_t^{-1}
\odot I_x(\theta_0)
\\
&\ge&\max_k B_t^{-1}\odot
I_k \ge\sum_k \alpha_k
\bigl(B_t^{-1}\odot I_k \bigr) =
B_t^{-1}\odot B^*
\\
&=& \tr \bigl(B_t^{-1}B^* \bigr) = \sum
_{i=1}^n \lambda_i = n + \sum
_{i=1}^n(\lambda_i-1) \ge n+\sum
_{i=1}^n\log(\lambda_i)
\\
&=& n+\log\det \bigl(B_t^{-1}B^* \bigr) = n+f \bigl(B^*
\bigr)-f(B_t),
\end{eqnarray*}
where $[\min\lambda_{B_t}\ge\mu]D_t\leadsto0$ for any $\mu>0$. Noting
that $\II_{X_{t+1}}(\theta_0) = \mathrm{E}_t(B\mid\theta_0)$, we obtain
\[
\mathrm{E}_t \bigl(f(B_{t+1})\mid\theta_0
\bigr)-f(B_t) \ge\frac{1}{t+1} \bigl( f \bigl(B^* \bigr) -
f(B_t) - D_{\mu,t} \bigr),
\]
where $D_{\mu,t} = D_t + (2nM^2\mu^{-2})/(t+1)$.

From now on, in order to keep the notation clean, we will implicitly
condition all probability statements on $\Theta=\theta_0$.

Let the constants $f_0<f_1<f(B^*)$ be arbitrary and define $\mu:=
\exp(f_0)M^{1-n}/2>0$. Suppose that some $t_0$ satisfies
$f(B_{t_0})\ge f_0$. Then, the definition of $\mu$ guarantees that
$\min\lambda_{B_{t_0}}\ge2\mu$. Let then
$\alpha\in\,]1,\exp(\mu/M) ]$ be arbitrary. Since
$\min\lambda_{B_{t}}$ can decrease by at most $M/t$ per each step,
we obtain
\[
\min\lambda_{B_t}\ge2\mu-\sum_{t=t_0+1}^{t_1}
\frac{M}t \ge2\mu- M\log\frac{t_1}{t_0} \ge\mu
\]
for all $t$ between $t_0$ and $t_1:= \lfloor\alpha t_0\rfloor$.
Thus, the following inequalities hold true for all
$t\in [t_0,t_1 [$:
%
%
\begin{eqnarray*}
\mathrm{E}_{t-1}f(B_t)-f(B_{t-1}) &\ge&
\frac{1}{t} \bigl( f \bigl(B^* \bigr) - f(B_{t-1}) -
D_{\mu,t-1} \bigr),
\\
\mathrm{E}_{t-1} \bigl(tf(B_t)-(t-1)f(B_{t-1})
\bigr) &\ge& f \bigl(B^* \bigr) - D_{\mu,t-1},
\\
\mathrm{E}_{t_0} \bigl(tf(B_t)-(t-1)f(B_{t-1})
\bigr) &\ge& f \bigl(B^* \bigr) - \mathrm {E}_{t_0}D_{\mu,t-1},
\\
\sum_{t=t_0+1}^{t_1}\mathrm{E}_{t_0}
\bigl(tf(B_t)-(t-1)f(B_{t-1}) \bigr) &\ge& \sum
_{t=t_0+1}^{t_1} \bigl(f \bigl(B^* \bigr) -
\mathrm{E}_{t_0}D_{\mu,t-1} \bigr),
\\
\mathrm{E}_{t_0} \bigl(t_1 f(B_{t_1}) \bigr) -
t_0 f(B_{t_0}) &\ge&(t_1-t_0)f
\bigl(B^* \bigr) - \sum_{t=t_0}^{t_1-1}
\mathrm{E}_{t_0}D_{\mu,t},
\end{eqnarray*}
and dividing by $t_1$, we obtain the inequality
%
%
\begin{eqnarray*}
\mathrm{E}_{t_0}f(B_{t_1}) - \alpha^{-1}f(B_{t_0})
&\ge& \biggl(1-\frac{t_0}{t_1} \biggr)f \bigl(B^* \bigr) -
\mathrm{E}_{t_0} \Biggl(\frac
{1}{t_1}\sum
_{t=t_0}^{t_1-1} D_{\mu,t} \Biggr)
\\
&\to& \bigl(1-\alpha^{-1} \bigr)f \bigl(B^* \bigr),
\end{eqnarray*}
where we have used the fact that $t_1\le\alpha t_0$, and where the
convergence holds for any increasing sequence of indices $t_0$
satisfying $f(B_{t_0})\ge f_0$ (which implies
$\min\lambda_{B_t}\ge\mu$ for all $t\in [t_0,t_1 [$). This
convergence is obtained by applying Lemma~\ref{lemproplemma}(3) to
the bounded sequence $[\min\lambda_{B_t}\ge\mu]D_{\mu,t}\leadsto0$,
which yields
\[
\Biggl\llvert \frac{1}{t_1}\sum_{t=t_0}^{t_1-1}
D_{\mu,t} \Biggr\rrvert \le \frac{1}{t_1}\sum
_{t=0}^{t_1-1} \bigl\llvert [\min\lambda_{B_t}
\ge\mu ]D_{\mu,t} \bigr\rrvert \to0
\]
(and since $\llvert  D_{\mu,t}\rrvert  \le2nM\mu^{-1} + 2nM^2\mu^{-2}$ for all $t$,
Lebesgue's dominated convergence theorem allows us to take this
limit inside the expectation). Thus, there exists a positive
constant $s$ such that
\[
\mathrm{E}_{t_0}f(B_{t_1}) \ge f(B_{t_0}) + 2s
\]
for all sufficiently large $t_0$ satisfying $f_0 \le f(B_{t_0})\le
f_1$. Also, since the maximum change in the value of $f$ over one
step is bounded by $v/t$ for some constant $v>0$ (depending on
$\mu$), we obtain
\[
\Var_{t_0}f(B_{t_1}) \le \sum_{t=t_0+1}^{t_1}
\biggl(\frac{v}t \biggr)^2 \le \int_{t_0}^{t_1}
\biggl(\frac{v}t \biggr)^2 \,\mathrm{d}t = v^2 \biggl(
\frac{1}{t_0}-\frac{1}{t_1} \biggr)\le\frac{v^2}{t_0}.
\]
Now Markov's inequality yields
%
%
\begin{eqnarray*}
\Pr_{t_0} \bigl\{f(B_{t_1}) < f(B_{t_0}) + s \bigr
\} &\le&\Pr_{t_0} \bigl\{f(B_{t_1}) < \mathrm{E}_{t_0}f(B_{t_1})
- s \bigr\}
\\
&\le&\Pr_{t_0} \bigl\{ \bigl\llvert \mathrm{E}_{t_0}f(B_{t_1})
- f(B_{t_1}) \bigr\rrvert ^2 > s^2 \bigr\}
\\
&\le&\frac{\Var_{t_0} f(B_{t_1})}{s^2} \le\frac{v^2}{t_0s^2}.
\end{eqnarray*}
As this upper bound on the probability sums to a finite number over
the sequence $t_0(k)$ determined by $t_0(k+1) = t_1(k)=\lfloor
\alpha t_0(k)\rfloor$, the Borel--Cantelli lemma implies that almost
surely $f(B_{t_0(k+1)}) < f(B_{t_0(k)}) + s$ holds for only finitely
many indices $k\in\mathbb{N}$ satisfying $f_0\le f(B_{t_0(k)})\le
f_1$. Thus, there exists $k_0$ such that for all $k\ge k_0$,
whenever $f_0\le f(B_{t_0(k)})\le f_1$, the value $f(B_{t_0}(k))$
will increase by at least $s$ on each step as $k$ increases.
Furthermore, since
\[
\bigl\llvert f(B_t)-f(B_{t_0(k)}) \bigr\rrvert \le\sum
_{t=t_0(k)+1}^{t_1(k)}\frac{v}t \le v\log
\frac{t_1(k)}{t_0(k)}\le v\log\alpha
\]
for all $t\in[t_0(k),t_1(k)[$, it follows that if $f(B_{t_0(k)})\ge
f_0$ for any $k\ge k_0$, then $f(B_t)\ge f_1-v\log\alpha$ for all
sufficiently large $t$ (provided that $f_1-v\log\alpha\ge f_0$).
Since $f_1-v\log\alpha$ can be made arbitrarily close to $f(B^*)$ by
appropriate choices of rational $\alpha>1$ and rational $f_1<f(B^*)$
for arbitrarily small rational $f_0$, we almost surely obtain
$\liminf_{t\to\infty} f(B_t)\ge f(B^*)$ unless $f(B_t)$ eventually
stays below any number. But this would imply that
$\limsup_{t\to\infty}\min\lambda_{B_t}\le0$, which is almost surely
contradicted by Lemma~\ref{lemminlambda}.
\end{pf}

%
\begin{corollary}[(Asymptotic D-optimality, part 2)]\label{cordopt}
Conditioned on almost any $\theta_0\in\mathsfTheta $ satisfying
\textup{O1--O4}, there exists a neighborhood $U$ of $\theta_0$ such that
$t\Cov_t(\Theta\mid U)\astoi(B^*)^{-1}$. This is optimal in the
sense that for any other strategy in place of \textup{O4} and \textup{C4}, almost
surely $\liminf_{t\to\infty}\det(t\Cov_t(\Theta\mid
U))\ge\det(B^*)^{-1}$.
\end{corollary}
\begin{pf}
Given O4, Theorems~\ref{thmdopt} and \ref{thmasymnorm}(2) imply
that $t\Cov_t(\Theta\mid U)\asto(B^*)^{-1}$. For any other
strategy, we have $\limsup_{t\to\infty} \det(B_t)\le\det(B^*)$ a.s.,
and so Theorem~\ref{thmasymnorm}(2) yields
$\liminf_{t\to\infty}\det(t\Cov_t(\Theta\mid U))\ge\det(B^*)^{-1}$
a.s. as $t$ increases within indices satisfying
$\min\lambda_{B_t}>\mu$ for some given $\mu>0$. But\vspace*{1pt}
Corollary~\ref{corvarlower} implies that if we choose a
sufficiently small $\mu>0$, then $\det(t\Cov_t(\Theta\mid
U))\ge\det(B^*)^{-1}$ also for $\min\lambda_{B_t}\le\mu$, and the
statement follows.
\end{pf}


\begin{remark}
As discussed in the beginning of this section, secondary modes with
weights proportional to $1/t$ may remain outside $U$, and they do
contribute to the asymptotic variance. Thus, the D-optimality
result (part 2) shown here is only a local form of optimality.

The situation would be different if the placements were chosen so as
to minimize the determinant of the posterior covariance
$\Cov_t(\Theta)$ directly (which, of course, presupposes that the
parameter space has global Euclidean structure). Then, slightly
more trials would be spent to decrease the weights of the secondary
modes, but they should remain insignificant in proportion. Thus, we
can conjecture that $B_t\asto B^*$ would still obtain in
Theorem~\ref{thmdopt} with $t\Cov_t(\Theta)$ asymptotically equal
to $(B_t)^{-1}$, making the result globally optimal.
\end{remark}

\subsection{Asymptotic entropy}\label{secasyment}

Here we use the D-optimality result to derive an expression for the
asymptotic entropy.

%
\begin{corollary}\label{cortinv}
Conditioned on almost any $\theta_0\in\mathsfTheta $ satisfying
\textup{O1--O4}, for any neighborhood $U$ of $\theta_0$, there exists a
constant $c_U$ such that almost surely, $p_t(U^c) \le c_U/t$ for
a.e. $t\in\mathbb{N}$.
\end{corollary}
\begin{pf}
Theorem~\ref{thmdopt} implies that $\min\lambda_{B_t}\ge\mu$ for
all sufficiently large $t$ for some $\mu>0$. Hence, given any
$\varepsilon>0$, Theorem~\ref{thmasymnorm}(3) yields
\[
t\I_t(\Theta;Y_{X_{t+1}}\mid U) \le \sup_{x\in\mathsf{X}}B_t^{-1}
\odot I_x(\theta_0) + \varepsilon\le n
\mu^{-1}M+\varepsilon=: c
\]
for all sufficiently large $t$, where $U$ is any sufficiently small
neighborhood of $\theta_0$. Combined with
Lemma~\ref{leminfolocal}, this implies that
$\I_t(\Theta;Y_{X_{t+1}})\le2c/t$ for a.e. $t\in\mathbb{N}$, and
so Lemma~\ref{lemscimpl}(2) yields the statement.
\end{pf}
%

\begin{remark}
Note that the statement of Corollary~\ref{cortinv} holds only for
a.e. $t\in\mathbb{N}$. What happens in a sufficiently long run is
that most trials are spent on increasing the accuracy around the
global mode and an approximately logarithmically growing number of
trials is spent on placements that decrease the weights of secondary
modes. However, on any such trial there is a small probability that
the weight of the secondary mode actually increases, and given a
sufficiently long run, this will eventually happen arbitrarily many
times in a row, making the weight of the secondary mode temporarily
arbitrarily much larger than the $c/t$ bound that holds on most
trials.
\end{remark}


\begin{theorem}\label{thmtinvp}
Conditioned on almost any $\theta_0\in\mathsfTheta $ satisfying
\textup{O1--O4}, if the prior entropy $\mathrm{H}(\Theta)$ w.r.t. a
parameterization that is consistent with the local Euclidean
structure (i.e., the prior density $p(\theta)$ is given w.r.t. a
measure that coincides with the Lebesgue measure on subsets of
$U_0$) is well-defined and finite, then, almost surely
\[
\mathrm{H}_t(\Theta) + \frac{n}2\log t \leadsto H^*:= -
\frac{1}2\log\det \bigl(B^* \bigr) + \frac{n}2\log(2\uppi \mathrm{e}).
\]
%
\end{theorem}

\begin{pf}
Let us condition everything on $\theta_0$ being the true value.
Theorem~\ref{thmasymnorm}(2) implies that for some sufficiently
small neighborhood $U$ of $\theta_0$,
\[
\mathrm{H}_t(\Theta\mid U) + \frac{n}2\log t \asto H^*.
\]
Lemmas \ref{lementkleq} and \ref{lementsublin} imply that for any
$\varepsilon>0$, $\llvert  \mathrm{H}_t(\Theta\mid U^c)\rrvert  <\varepsilon t$ for
all sufficiently
large $t$, and as Corollary~\ref{cortinv} yields $p_t(U^c)\le c/t$
for a.e. $t$, Lemma~\ref{lemproplemma}(2) implies
$p_t(U^c)\mathrm{H}_t(\Theta\mid U^c)\leadsto0$. The statement now follows
from the chain rule of entropy
\[
\mathrm{H}_t(\Theta) = p_t(U)\mathrm{H}_t(
\Theta\mid U) + \underbrace{p_t \bigl(U^c \bigr)
\mathrm{H}_t \bigl(\Theta\mid U^c \bigr)}_{\leadsto0}+
\underbrace{\mathrm{H}_t \bigl([\Theta\in U] \bigr)}_{\to0~\mathrm
{a.s.}},
\]
where the first term satisfies
\[
p_t(U)\mathrm{H}_t(\Theta\mid U) + \frac{n}2
\log t = p_t(U) \biggl[\mathrm{H}_t(\Theta\mid U) +
\frac{n}2\log t \biggr] +\underbrace{p_t
\bigl(U^c \bigr)}_{\le c/t}\frac{n}2\log t\leadsto H^*.
\]
\upqed
\end{pf}

%
\begin{corollary}\label{corasyment}
Suppose that \textup{O1--O4} hold for almost all $\theta_0\in\mathsfTheta $
and that the prior entropy $\mathrm{H}(\Theta)$ w.r.t. a parameterization
that is consistent with the local Euclidean structures $U_0$ in \textup{O2}
is well-defined and finite. Then,
\[
\mathrm{H}_t(\Theta) + \frac{n}2\log t \prhoto H^*.
\]
In other words, there exists a set $K\subset\mathbb{N}$ of indices
with $\rho(K)=1$ such that
\[
\mathrm{H}_t(\Theta) + \frac{n}2\log t \pto H^*,
\]
as $t$ increases within $K$.
\end{corollary}
\begin{pf}
Apply Lemma~\ref{lemasprhoto} to the statement of Theorem~\ref{thmtinvp}.
\end{pf}

\subsection{Varying cost of observation}\label{secvaryingcost}

In Kujala \cite{bestvalue} the adaptive sequential estimation framework is
generalized to the situation where the observation of $Y_x$ is
associated with some random cost $C_x$ of observation, which given the
value of $Y_x$, is independent of $\Theta$ and the results and costs
of any other observations:
\par
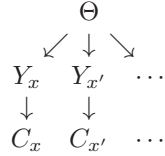
\[
\begin{tabular} {c@{}c@{}c@{}c@{}c} & & $\Theta$
\\
& $\swarrow$ & $\downarrow$ & $\searrow$
\\
$Y_x$ & & $Y_{x'}$ & & $\cdots$
\\
$\downarrow$ & & $\downarrow$
\\
$C_x$ & & $C_{x'}$ & & $\cdots$ \end{tabular} %
\]
The technical requirement that $C_x$ depends on $\Theta$ only through
$Y_x$ is satisfied in particular if $C_x$ is a component of $Y_x$.
Thus, it leads to no loss of generality if the incurred costs are
observable.

The goal considered in Kujala \cite{bestvalue} is maximization of the
expected information gain of a sequential experiment that terminates
when the total cost overruns a given budget. To achieve this goal,
the heuristic of maximizing the expected information gain
$\I_t(\Theta;Y_x)$ divided by the expected cost $\mathrm{E}_t(C_x)$
on each
trial is proposed. In this section, we are able to show that this
heuristic is in fact asymptotically optimal (as the budget tends to
infinity) under essentially the same conditions that the plain
information gain maximization is.

Thus, condition O4 is now replaced by the following:
\begin{enumerate}[O4$'$.]
\item[O4$'$.] The placements satisfy
\[
R'_t:=\frac{
{\I_t(\Theta;Y_{X_{t+1}})}/{\mathrm{E}_t(C_{X_{t+1}})}
}{
\sup_{x\in\mathsf{X}}({\I_t(\Theta;Y_x)}/{\mathrm{E}_t(C_x)})
}\leadsto1,
\]
where $\llvert  C_x\rrvert  \le M$, $\mathrm{E}(C_x\mid\theta_0)\ge\gamma'>0$, and the
family of expected cost functions $\{ \theta\mapsto
\mathrm{E}(C_x\mid\theta)\dvt x\in\mathsf{X} \}$ is equicontinuous at
$\theta_0$.
\end{enumerate}
Due to the assumed bounds on the expected cost $\mathrm{E}(C_x\mid
\theta_0)$,
condition C4 is still satisfied and so all the previous lemmas
depending on it apply. Together with the following lemma, these
bounds also imply that the total cost grows asymptotically within
linear bounds.

\begin{lemma}\label{lemasymcost}
Suppose that \textup{O4$'$} holds. Then, conditioned on $\theta_0$ as the
true parameter value,
\[
\frac{C_t -\sum_{k=1}^t \mathrm{E}(C_{X_k}\mid\theta_0)}{t}\astoi0,
\]
where $C_t:= \sum_{k=1}^t C_{X_k}$. In particular, for any
$\gamma<\gamma'$, almost surely $C_t\ge t\gamma$ for all
sufficiently large $t$ (as well as $C_t\le tM$ for all $t$).
\end{lemma}

\begin{pf}
Denoting $Z_k = C_{X_k} - \mathrm{E}(C_{X_k}\mid\theta_0)$, given
$\Theta=\theta_0$, the sequence $Z_1+\cdots+Z_k$ of partial sums is a
martingale and satisfies $\mathrm{E}(\llvert  Z_k\rrvert  ^2)\le M^2 <\infty$ for all $k$,
and so Theorem~\ref{thmmartingaleSLLN} implies that
$(Z_1+\cdots+Z_t)/t\asto0$, which is the statement.
\end{pf}
Next, we will generalize Corollary~\ref{corinfohess} for the
cost-aware placements.

\begin{corollary}\label{corinfohess2}
Conditioned on almost any $\theta_0$ satisfying \textup{O1--O3} and \textup{O4$'$},
the sequence
\[
D_t:= \sup_{x\in\mathsf{X}}B_t^{-1}
\odot\frac{I_x(\theta
_0)}{\mathrm{E}(C_x\mid\theta_0)} - B_t^{-1}\odot\frac{\II_{X_{t+1}}(\theta_0)}{\mathrm
{E}_t(C_{X_{t+1}}\mid\theta_0)}
\]
satisfies $[\min\lambda_{(C_t/t)B_t}\ge\mu]D_t\leadsto0$ a.s. for
any given $\mu>0$, where $\min\lambda_{(C_t/t)B_t}$ denotes the
smallest eigenvalue of $B_t:= -C_t^{-1}\nabla^2_\theta \log
p(\mathbf{Y}_t\mid\theta_0)$ and $C_t:= \sum_{k=1}^t C_{X_k}$.
\end{corollary}

\begin{pf}
Let us first shrink the neighborhood $U_0$ of $\theta_0$ as
necessary to make its diameter smaller than the constant
$\delta_{\mu,C}$ given by Lemma~\ref{lemL12}. Then, let $U\subset
U_0$ be the neighborhood of $\theta_0$ given by
Lemma~\ref{leminfolocal}. The boundedness and equicontinuity at
$\theta_0$ of $\theta\mapsto\mathrm{E}(C_x\mid\theta)\in[\gamma
',M]$ imply
that conditioned on $\Theta=\theta_0$, almost surely, $\mathrm
{E}_t(C_x)\to
\mathrm{E}(C_x\mid\theta_0)$, uniformly over all $x\in\mathsf{X}$. Combined
with Theorem~\ref{thmasymnorm}(3), this implies that there exist
random sequences $E_t\to0$ and $E_t'\to0$ such that conditioned on
$\theta_0$ as the true value,
%
%
\begin{eqnarray*}
\frac{1}2 \sup_{x\in\mathsf{X}} B_t^{-1}
\odot\frac{I_x(\theta
_0)}{\mathrm{E}(C_x\mid\theta_0)} &=& \sup_{x\in\mathsf{X}} C_t
\frac{\I_t(\Theta;Y_x\mid
U)}{{\mathrm{E}_t(C_x)}} + E_t,
\\
\frac{1}2 B_t^{-1}\odot\frac{\II_{X_{t+1}}(\theta_0)}{\mathrm{E}(C_{X_{t+1}}\mid\theta_0)} &=&
C_t\frac{\I_t(\Theta;Y_{X_{t+1}}\mid U)}{{\mathrm{E}_t(C_{X_{t+1}})}} + E_t'
\end{eqnarray*}
whenever $\min\lambda_{(C_t/t)B_t}\ge\mu$. For these $t$, it follows
%
%
\begin{eqnarray*}
\frac{1}2D_t &=& \biggl(\frac{1}2 \sup
_{x\in\mathsf{X}} \underbrace{B_t^{-1}\odot
\frac{I_x(\theta_0)}{\mathrm{E}(C_x\mid\theta_0)}}_{
\le\tr(B_t^{-1}I_x(\theta_0))/\gamma\le n(\gamma\mu)^{-1}M
} - E_t \biggr) \biggl(1 -
\frac{
{\I_t(\Theta;Y_{X_{t+1}}\mid U)}/{\mathrm{E}_t(C_{X_{t+1}})}
}{
\sup_{x\in\mathsf{X}} ({\I_t(\Theta;Y_x\mid U)}/{\mathrm{E}_t(C_x)})
} \biggr)
\\
&&{} + E_t-E_t',
\end{eqnarray*}
where Lemma~\ref{leminfolocal} and the inequality
$\I_t(\Theta;Y_x)\ge p_t(U)\I_t(\Theta;Y_x\mid U)$ yield
%
%
\begin{eqnarray*}
\frac{
{\I_t(\Theta;Y_{X_{t+1}}\mid U)}/{\mathrm{E}_t(C_{X_{t+1}})}
}{
\sup_{x\in\mathsf{X}} ({\I_t(\Theta;Y_x\mid U)}/{\mathrm{E}_t(C_x)})
} \ge p_t(U)\frac{
{\I_t(\Theta;Y_{X_{t+1}}\mid U)}/{\mathrm{E}_t(C_{X_{t+1}})}
}{
\sup_{x\in\mathsf{X}} ({\I_t(\Theta;Y_x)}/{\mathrm{E}_t(C_x)})
} = p_t(U)
Q_t R'_t \leadsto1,
\end{eqnarray*}
and so $[\min\lambda_{(C_t/t)B_t}\ge\mu]D_t\leadsto0$.
\end{pf}

%
\begin{lemma}\label{leminfocostrange}
The range of the expression
\[
r_t = \frac{\sum_{k=1}^t I_{x_k}(\theta_0)}{\sum_{k=1}^t \mathrm
{E}(C_{x_k}\mid\theta_0)}
\]
over all sequences $x_k$ in $\mathsf{X}$ and all finite $t$ is a
dense subset of the set $\mathcal{I}$ defined as the closure of the
convex hull of
\[
S = \biggl\{\frac{I_x(\theta_0)}{\mathrm{E}(C_x\mid\theta
_0)} \biggr\}_{x\in\mathsf{X}}.
\]
Furthermore, the range of the limits of all converging $r_t$ equals
$\mathcal{I}$.
\end{lemma}
\begin{pf}
For any sequence $\{x_k\}$, we have
\[
r_t = \frac{\sum_{k=1}^t I_{x_k}(\theta_0)}{\sum_{k=1}^t \mathrm
{E}(C_{x_k}\mid\theta_0)} = \sum_{k=1}^t
\underbrace{ \biggl(\frac{\mathrm{E}(C_{x_k}\mid\theta
_0)}{{\sum_{k=1}^t \mathrm{E}(C_{x_k}\mid\theta_0)}} \biggr)}_{=:\alpha_{k,t}}
\frac{I_{x_k}(\theta_0)}{\mathrm{E}(C_{x_k}\mid
\theta_0)},
\]
and so $r_t$ is always a convex combination of elements in $S$. The
convex combination is not exactly linear w.r.t. the number of
different $x$ in the sequence because of the different
$\mathrm{E}(C_{x_k}\mid\theta_0)$ weights, but nonetheless, by
varying the
proportions of different $x$ in a sufficiently long sequence, any
convex combination can be approximated arbitrarily well.
\end{pf}


\begin{theorem}[(Asymptotic D-optimality, part 1)]\label{thmdopt-myopic}
Conditioned on almost any $\theta_0\in\mathsfTheta $ satisfying
\textup{O1--O3}, \textup{O4$'$}, almost surely,
%
%
\begin{eqnarray*}
B_t:= \frac{-\nabla^2_\theta \log p(\mathbf{Y}_t\mid\theta_0)}{C_t} \to B^*:= \argmax_{B\in\mathcal{I}}
\det(B),
\end{eqnarray*}
where $C_t:= \sum_{k=1}^tC_{X_{k}}$ and $\mathcal{I}$ is the convex
hull of the closure of
\[
S = \biggl\{ \frac{I_x(\theta_0)}{\mathrm{E}(C_x\mid\theta
_0)}\dvt x\in\mathsf{X} \biggr\}.
\]
This is optimal in the sense that for any strategy of choosing the
placements $X_t$ (instead of \textup{O4$'$} and \textup{C4}), almost surely
$\limsup_{t\to\infty}\det(B_t)\le\det(B^*)$.
\end{theorem}

\begin{pf}
Since $S$ is bounded, $\mathcal{I}$ is a compact convex set and
$B^*$ is well defined. Lemmas \ref{lemasymhess},
\ref{lemasymcost}, and \ref{leminfocostrange} imply that
$\limsup_{t\to\infty}\det(B_t)\le\det(B^*)$ a.s. Let us then show
that this upper bound is tight.

Lemma \ref{leminfocostrange} implies that there exists a representation
\[
B^* = \lim_{m\to\infty}\frac{\sum_{k=1}^mI_k}{\sum_{k=1}^mc_k}
\]
of the optimum point $B^*$ where $(I_k,c_k)$ are elements of
$\{ (I_x(\theta_0),\mathrm{E}(C_x\mid\theta_0))\dvt x\in\mathsf{X}
\}$.

Denoting $B:= -\nabla^2_\theta \log(p(Y_{X_{t+1}}\mid\theta_0))$
and $C:= C_{X_{t+1}}$, and assuming
$\min\lambda_{(C_t/t)B_t}\ge\mu$, we obtain
\[
\llvert B\rrvert,\llvert C\rrvert \le M,\qquad \bigl\llvert B_t^{-1}
\bigr\rrvert \le(\mu/M)^{-1}, \qquad \llvert B-CB_t\rrvert
\le M+M^2/\mu,\qquad C_t+C\ge\gamma(t+1)
\]
and so, for some $B'$ between $0$ and $B_{t+1}-B_t$, we obtain
%
%
\begin{eqnarray*}
f(B_{t+1}) - f(B_t) &=& f \biggl(\frac{C_tB_t+B}{C_t+C} \biggr)
- f(B_t)
\\
&=& B_t^{-1}\odot\frac{B-CB_t}{C_t+C} -\frac{1}2\tr
\bigl(B_t^{-1}B'B_t^{-1}B'
\bigr)
\\
&\ge&\frac{1}{C_t+C} \biggl(B_t^{-1}\odot B - nC -
\frac{[(\mu
/M)^{-1}(M+M^2/\mu)]^2}{C_t+C} \biggr)
\\
&\ge&\underbrace{\frac{\mathrm{E}_t(C\mid\theta_0)}{C_t+C}}_{\ge
(\gamma/M)/(t+1)} \biggl(B_t^{-1}
\odot\frac{B}{\mathrm{E}_t(C\mid
\theta_0)} - \frac{nC}{\mathrm{E}_t(C\mid\theta_0)} - \frac
{C_{M,\mu,\gamma}}{t+1} \biggr).
\end{eqnarray*}
Denoting by $\lambda_i$ the eigenvalues of $B_t^{-1}B^*$, we obtain
%
%
\begin{eqnarray*}
B_t^{-1}\odot\frac{\I_{X_{t+1}}(\theta_0)}{\mathrm{E}(
C_{X_{t+1}}\mid\theta_0)} + D_t &=& \sup
_{x\in\mathsf{X}}B_t^{-1}\odot\frac{I_x(\theta_0)}{\mathrm
{E}(C_x\mid\theta_0)}
\\
&\ge&\sup_k \biggl(B_t^{-1}\odot
\frac{I_k}{c_k} \biggr) \ge\lim_{m\to\infty}\frac{\sum_{k=1}^m(B_t^{-1}\odot I_k)}{\sum_{k=1}^mc_k} =
B_t^{-1}\odot B^*
\\
&=& \tr \bigl(B_t^{-1}B^* \bigr) = \sum
_{i=1}^n \lambda_i = n + \sum
_{i=1}^n(\lambda_i-1) \ge n+\sum
_{i=1}^n\log(\lambda_i)
\\
&=& n+\log\det \bigl(B_t^{-1}B^* \bigr) = n+f \bigl(B^*
\bigr)-f(B_t),
\end{eqnarray*}
where Corollary~\ref{corinfohess2} implies that
$[\min\lambda_{(C_t/t)B_t}\ge\mu]D_t\leadsto0$. Noting that
$\mathrm{E}_t(B/\mathrm{E}_t(C\mid\theta_0)\mid\theta_0) =
\I_{X_{t+1}}(\theta_0)/\mathrm{E}( C_{X_{t+1}}\mid\theta_0)$, it follows
\[
\mathrm{E}_t \bigl(f(B_{t+1})\mid\theta_0
\bigr) - f(B_t) \ge\frac{\gamma
/M}{t+1} \bigl(f \bigl(B^*
\bigr)-f(B_t)-D_{\mu,t} \bigr),
\]
where $D_{\mu,t} = D_t + C_{M,\mu,\gamma}/(t+1)$.

From here on, the proof is essentially the same as in the maximum
information case. We just use $\mu:= \exp(f_0)M^{-n}/2$ to
guarantee that $\min\lambda_{(C_t/t)B_t}\ge2\mu$ for $f(B_t)\ge
f_0$.
\end{pf}

The part 2 of the D-optimality result as well as analogs of the
asymptotic entropy results follow with essentially the same proofs
(just replacing $t$ with $C_t$ at appropriate places):


\begin{corollary}[(Asymptotic D-optimality, part 2)]
Conditioned on almost any $\theta_0\in\mathsfTheta $ satisfying
\textup{O1--O3}, \textup{O4$'$}, there exists a neighborhood $U$ of $\theta_0$ such
that $C_t\Cov_t(\Theta\mid U)\astoi(B^*)^{-1}$, where $C_t:=
\sum_{k=1}^t C_{X_k}$. This is optimal in the sense that for any
other strategy in place of \textup{O4$'$} and \textup{C4}, almost surely
$\liminf_{t\to\infty}\det(C_t\Cov_t(\Theta\mid
U))\ge\det(B^*)^{-1}$.
\end{corollary}


\begin{theorem}
Conditioned on almost any $\theta_0\in\mathsfTheta $ satisfying
\textup{O1--O3}, \textup{O4$'$}, if the prior entropy $\mathrm{H}(\Theta)$ w.r.t. a
parameterization that is consistent with the local Euclidean
structure (i.e., the prior density $p(\theta)$ is given w.r.t. a
measure that coincides with the Lebesgue measure on subsets of
$U_0$) is well-defined and finite, then, almost surely
\[
\mathrm{H}_t(\Theta) + \frac{n}2\log C_t
\leadsto H^*:= -\frac{1}2\log\det \bigl(B^* \bigr) + \frac{n}2
\log(2\uppi \mathrm{e}),
\]
where $C_t:= \sum_{k=1}^t C_{X_k}$.
\end{theorem}

%
\begin{corollary}\label{corasyment-myopic}
Suppose that \textup{O1--O4} hold for almost all $\theta_0\in\mathsfTheta $
and that the prior entropy $\mathrm{H}(\Theta)$ w.r.t. a parameterization
that is consistent with the local Euclidean structures $U_0$ in \textup{O2}
is well-defined and finite. Then,
\[
\mathrm{H}_t(\Theta) + \frac{n}2\log C_t
\prhoto H^*,
\]
where $C_t:= \sum_{k=1}^t C_{X_k}$. In other words, there exists a
set $K\subset\mathbb{N}$ of indices with $\rho(K)=1$ such that
\[
\mathrm{H}_t(\Theta) + \frac{n}2\log C_t \pto
H^*,
\]
as $t$ increases within $K$.
\end{corollary}

\section{Examples}\label{secexamples}

In this section, we give specific examples illustrating the optimality
results.

%
\begin{example}[(Psychometric model)]
Consider the psychometric model, where an observer's unknown intensity
threshold $\Theta$ for detecting a stimulus of intensity $x$ is
distributed uniformly on $[0,100]$ and the trial result $Y_x\in\{0,1\}$
for a test intensity $x\in[0,100]$ is distributed as
\[
p(y_x\mid\theta)=\cases{\displaystyle \psi(\theta-x), &\quad
$y_x=1$ (detected),
\cr
\displaystyle 1-\psi(\theta-x), &\quad
$y_x=0$ (not detected),}
\]
where $\psi(x)$ is the psychometric function, here assumed to be the sigmoid
\[
\psi(x)=\frac{1}{1+\mathrm{e}^{-x}}
\]
for simplicity (for more general psychometric models, see
Kujala and Lukka \cite{kujalalukka2006}, and the references therein).

In this model, the Fisher information of a given placement $x$ is
calculated as
\[
I_x(\theta)=\sum_{y_x=0}^1p(y_x
\mid\theta) \biggl[\frac{\partial
}{\partial\theta}\log p(y_x\mid\theta)
\biggr]^2 =\frac{\psi'(\theta-x)^2}{\psi(\theta-x)[1-\psi(\theta-x)]} =\frac{\mathrm{e}^{\theta-x}}{ [1+\mathrm{e}^{\theta-x} ]^2}.
\]
Thus, for any given $\theta_0$, the D-optimal value of the averaged
Fisher information in Theorem~\ref{thmdopt} is $B^*=\frac{1}{4}$
given by the placement $x=\theta_0$ to which the greedy algorithm
eventually converges. Now Corollary~\ref{corasyment} yields
%
%
\begin{equation}
\label{eqpsient-optimal} \mathrm{H}_t(\Theta)+\frac{n}{2}\log t\prhoto
H^* =-\frac{1}{2}\log\underbrace{\det \bigl(B^* \bigr)}_{=0.25}+\,
\frac{n}{2}\log (2\uppi \mathrm{e})
\end{equation}
and this is the asymptotically optimal posterior entropy. In this
example, the same expression also gives the asymptotically optimal
expected utility $\mathrm{E}(\mathrm{H}_t(\Theta))+\frac{n}2\log t$,
which we
will next compare to that of the offline design.
\end{example}

%
\begin{example}[(Offline design)]\label{exuniform-placement}
A rigorous study of the optimal offline design is beyond the scope of
the present article, so we will not go into detailed proofs here but
only sketch the general ideas. Suffice it to say that for an offline
design for optimizing the expected utility $\mathrm{E}(\mathrm
{H}_t(\Theta))$, one
cannot do much better than to use the usual strategy of placing the
trials evenly on the interval $[0,100]$. (Due to boundary effects, an
exactly uniform distribution of placements is not really the global
optimum, but for simplicity, we avoid a more complicated discussion
here.)

For uniform placement of trials on $[0,100]$, Lemma~\ref{lemasymhess}
implies
\[
B_t\asto\frac{1}{100}\int_0^{100}I_x(
\theta_0)\,\mathrm{d}x =\frac{1}{100} \biggl(\frac{1}{1+\mathrm{e}^{-\theta_0}}-
\frac
{1}{1+\mathrm{e}^{100-\theta_0}} \biggr) \in[0.005,0.01],
\]
where $B_t=-t^{-1}\nabla_{\theta}^2\log p(\mathbf{Y}_t\mid\theta_0)$,
and it can be shown that the asymptotic posterior entropy satisfies
\[
\mathrm{H}_t(\Theta)+\frac{n}{2}\log t - \biggl[-
\frac{1}{2}\log\underbrace{\det(B_t)}_{\lim\le0.01} \biggr] +
\frac{n}{2}\log(2\uppi \mathrm{e})\asto0,
\]
which implies the asymptotic lower bound
\[
\liminf_{t\to\infty} \biggl[ \mathrm{H}_t(\Theta)+
\frac{n}{2}\log t \biggr] \ge-\frac{1}2\log0.01 +
\frac{n}2\log(2\uppi \mathrm{e})
\]
on\vspace*{1pt} the posterior entropy.
Comparing to the asymptotically optimal posterior entropy
(\ref{eqpsient-optimal}), it follows that the offline design needs
asymptotically at least $(\frac{0.25}{0.01})^{{1}/{n}}=25$
times as many trials as the optimal adaptive design for the same
accuracy. If the range $[0,100]$ is doubled, then this number
approximately doubles as well, so the gap to the asymptotically
optimal adaptive design can be arbitrarily large.
\end{example}

%
\begin{example}[(Varying cost of observation)]
Let us then return to the adaptive case and suppose that instead of a
unit cost, each trial costs
\[
C_x=1+3[Y_x=0]
\]
units. Such a formulation could be based on the assumption that the
observer takes four times as long to respond when the stimulus is not
detected. Then, the asymptotic efficiency of a placement~$x$ in
Theorem~\ref{thmdopt-myopic} is characterized by the expression
%
%
\begin{equation}
\label{eqexample-optimal-cost-aware} \frac{I_x(\theta_0)}{\mathrm{E}(C_x)}=\frac{I_x(\theta
_0)}{1+3[1-\psi(\theta_0-x)]}=\frac{1}{5+5\cosh(\theta_0-x)-3\sinh
(\theta_0-x)}.
\end{equation}
This expression is maximized by the placement $x=\theta_0+\log2$ to
which the myopic algorithm eventually converges to (provided it is
within the range $[0,100]$). Thus, assuming that
$\theta_0\le100-\log2\approx99.3069$ and substituting the maximizer in
(\ref{eqexample-optimal-cost-aware}), we obtain in
Theorem~\ref{thmdopt-myopic} the D-optimal asymptotic efficiency
$B^*=\frac{1}{9}$. Comparing to the asymptotically optimal
placement $x=\theta_0$ for unit cost (yielding $B^*=\frac{1}{10}$
in (\ref{eqexample-optimal-cost-aware})), we see that the cost-aware
strategy reaches the same accuracy in $10\%$ less cost (time) in this
example.
\end{example}

\section{Discussion}\label{secdiscussion}

We have derived an expression for the asymptotic efficiency of any
sequential experiment design for both the standard framework with unit
cost of observation as well as for the generalized framework with
random costs of observation as proposed in Kujala \cite{bestvalue}. We have
shown an asymptotic D-optimality result for the greedy information
optimization strategy in the standard framework and we have extended
this result for the novel myopic strategy proposed in Kujala \cite{bestvalue}
for the situation with random costs of observations. These results
indicate that for (almost) all true parameter values $\theta_0$, the
greedy or myopic adaptive design is asymptotically optimal among all
placement strategies in a well-defined sense.

Assuming the standard sequential estimation framework with unit cost
of observation, Lemma~\ref{lemasymhess} together with the asymptotic
normality result imply that the asymptotic efficiency of any given
design is characterized by the average
\[
\frac{\sum_{k=1}^t I_{X_k}(\theta_0)}{t}
\]
of the Fisher information matrices $I_x(\theta_0)$ over the sequence
of placements $X_t$ and the \mbox{D-}optimality criterion of a design refers
to maximality of the determinant of this averaged information matrix
at the limit. For any given $\theta_0$, there is a distribution (or
sequence) of placements $x\in\mathsf{X}$ yielding the D-optimal
average information matrix. For (almost) all $\theta_0$, the
placements of the greedy adaptive design converge to such an optimum,
whereas the offline design cannot adjust the distribution of the
placements $x\in\mathsf{X}$ depending on the true value $\theta_0$.
Thus, the offline design can be equally efficient for a given true
value of $\Theta$, but generally not for all values
$\theta_0\in\mathsfTheta $ and depending on the model, the gap in
efficiency can be arbitrarily large as seen in
Example~\ref{exuniform-placement}.

The situation is essentially the same in the framework with random
costs of observation, the only difference being that the convergence
of the estimate of $\Theta$ is not measured in relation to $t$ but in
relation to the total cost $C_t = C_{X_1}+\cdots+C_{X_t}$ of
placements. In this situation, the asymptotic efficiency is
characterized by the ratio
\[
\frac{\sum_{k=1}^t I_{X_k}(\theta_0)}{\sum_{k=1}^t\mathrm
{E}(C_{X_k}\mid\theta_0)}
\]
and the limit is again determined by the distribution (or sequence) of
the placements $x\in\mathsf{X}$. Theorem~\ref{thmdopt-myopic} shows
that the myopic strategy of maximizing
\[
\frac{\I_t(\Theta;Y_x)}{\mathrm{E}_t(C_x)}
\]
yields the asymptotically D-optimal efficiency in this situation.

However, the actual utility function assumed in both of the frameworks
considered is the differential entropy, and so the most relevant
asymptotic optimality criterion should be based on the asymptotic
properties of the differential entropy as shown in, for example,
Corollaries \ref{corasyment}~and~\ref{corasyment-myopic}. Thus, a
topic for future work is finding conditions under which the results of
Corollaries~\ref{corasyment} and \ref{corasyment-myopic} can be
said to be optimal among all placement strategies.

\begin{appendix}
\section*{Appendix: Auxiliary theorems}

%
\begin{theorem}[(Stone--{\v C}ech compactification)] \label{thmstonecech}
Suppose that $X$ is a Tychonoff space. Then there exists a compact
space $\beta X$ that embeds $X$ as a dense subspace. Any continuous
map $f\dvtx X\to K$, where $K$ is a compact Hausdorff space, lifts
uniquely to a continuous map $\beta X\to K$.
\end{theorem}

%
\begin{theorem}[(Martingale convergence)]\label{thmmartingaleconv}
Let $X_k$ be a submartingale (i.e., $\mathrm{E}(X_{k+1}\mid
X_1,\ldots, X_k)\ge X_k$) and suppose that $\sup_k\mathrm
{E}\llvert  X_k\rrvert  <\infty$.
Then, $X = \lim_{k\to\infty}$ exists almost surely and
$\mathrm{E}\llvert  X\rrvert  <\infty$.
\end{theorem}

\begin{pf}
For example, \cite{schervish1995}, Theorem~B.117, page~648,  or
\cite{shiryaev1996}, Theorem~1, page~508.
\end{pf}


%
\begin{theorem}[(A strong law of large numbers for martingales)]\label{thmmartingaleSLLN}
Let $X_k = Z_1+\cdots+Z_k$ be a martingale and let $\delta>0$. If
\[
\sum_{k=1}^\infty\frac{\mathrm{E}(\llvert  Z_k\rrvert  ^2)}{k^{2\delta}}<\infty,
\]
then $X_k/k^\delta\astoi0$.
\end{theorem}
\begin{pf}
For example, \cite{chow1967} or \cite{shiryaev1996}, Theorem~4, page~519.
\end{pf}


\begin{theorem}[(Hoeffding--Azuma inequality)]\label{thmazuma}
Let $X_k$ be a martingale and suppose that $\llvert  X_k-X_{k-1}\rrvert  \le c_k$
for all $k$. Then, for all $t>0$ and $k\in\mathbb{N}$,
\[
\Pr\{X_n-X_0\ge t\}\le\exp \biggl(-\frac{t^2}{2\sum_{k=1}^nc_k^2}
\biggr),
\]
and
\[
\Pr \bigl\{\llvert X_n-X_0\rrvert \ge t \bigr\}\le2\exp
\biggl(- \frac{t^2}{2\sum_{k=1}^nc_k^2} \biggr).
\]
\end{theorem}

\begin{pf}
See \cite{hoeffding1963}, Theorem 2 and note around (2.18) on page~18,
or \cite{azuma1967}.
\end{pf}
\end{appendix}

\section*{Acknowledgements}
The author is grateful to Matti Vihola for many stimulating
discussions. This research was supported by the Academy of Finland
(grant number 121855).





%

\printhistory
\end{document}